\documentclass[a4paper, 11pt]{article}
\usepackage{fullpage,amsmath, amssymb,amscd}
\usepackage[mathscr]{eucal}

\newcommand\cyr{%
 \renewcommand\rmdefault{wncyr}%
 \renewcommand\sfdefault{wncyss}%
 \renewcommand\encodingdefault{OT2}%
\normalfont\selectfont} \DeclareTextFontCommand{\textcyr}{\cyr}

\newtheorem{theorem}{Theorem}
\newtheorem{lemma}[theorem]{Lemma}

\newtheorem{proposition}[theorem]{Proposition}

\def\mod{\operatorname{mod}}
\def\O{\operatorname{O}}

\def\GL{\operatorname{GL}}
\def\tr{\operatorname{tr}}
\def\Aut{\operatorname{Aut}}

\def\ds{\displaystyle}

\newenvironment{proof}{\par\noindent{\bf Proof.}}{$\square$\par\bigskip}

\newcommand{\const}{{\mathfrak{C}}}

\newcommand{\qed}{{$\square$\par\bigskip}}

\newcommand{\bbz}{{\mathbb{Z}}}

\newcommand{\F}{{\mathbb{F}}}
\newcommand{\Q}{{\mathbb{Q}}}
\newcommand{\Z}{{\mathbb{Z}}}

\newcommand{\leg}[2]{ \left( \frac{#1}{#2} \right) }

\begin{document}

\title{Average twin prime conjecture for elliptic curves}

\author{
Antal Balog (Alfr\'ed R\'enyi Institute of Mathematics)
\\ {\small{ balog@renyi.hu }}
\\
Alina-Carmen  Cojocaru  (University of Illinois at Chicago)
\thanks{Research supported in part by NSF
grant DMS-0636750.}
\\
{\small{ cojocaru@math.uic.edu }}
\\
Chantal David (Concordia University) \thanks{Research supported in
part by NSERC.}\\ {\small{ cdavid@mathstat.concordia.ca }} }

 \maketitle

\begin{abstract}
Let $E$ be an elliptic curve over $\Q$. In 1988,  Koblitz conjectured a
precise asymptotic for the number of primes $p$ up to $x$ such that
the order of the group of points of $E$ over $\F_p$ is prime. This
is an analogue of the Hardy and Littlewood twin prime conjecture in
the case of elliptic curves.

Koblitz's conjecture is still widely open.
In this paper we prove that Koblitz's conjecture is true on average
over a two-parameter family of elliptic curves. One of the key
ingredients in the proof is  a short average distribution result in the style  of
Barban-Davenport-Halberstam,  where the average is taken  over twin primes and their differences.
\end{abstract}

\tableofcontents

\renewcommand\baselinestretch{1.1}

\section{Introduction}

\hspace*{0.5cm} A well-known open problem in number theory is
{\it{the twin prime conjecture}}, which states that
that there exist infinitely many primes $p$ such that
$p+2$ is also a prime.
This conjecture was generalized by Alphonse de Polignac in 1849 to the statement that,
 for any even integer $r \neq
0$, there exist infinitely many primes $p$ such that $p+r$ is also a
prime. In 1922,  G.H. Hardy and J. Littlewood made this statement precise,
 predicting that, as $x \rightarrow \infty$,
$$
\#\{p \leq x: p+r \; \text{is prime}\}
\sim
{\mathfrak S}(r) \frac{x}{\log^2 x },
$$
where
\begin{equation}
\label{twin-prime-constant}
{\mathfrak S}(r)
:=
\left\{ \begin{array}{cc} \displaystyle{
2 \ds\prod_{\ell \neq 2} \frac{\ell(\ell-2)}{(\ell-1)^2} \prod_{\ell \mid r, \ell \neq 2}
\frac{\ell-1}{\ell-2}}
& \mbox{if $2 \mid r$},
\\
0 & \mbox{otherwise}.
\end{array} \right.
\end{equation}
Here and everywhere in the paper, $p$ and $\ell$ are used to denote primes.

Even though still inaccessible by current methods, the twin prime
conjecture has generated tremendous advances in number theory.
Indeed, in 1919 Viggo Brun \cite{Br} developed what is now known as the Brun
sieve to prove the surprising result that
$\ds\sum_{p \atop{p+2 \; \text{prime}}} \frac{1}{p} < \infty$.
Brun's methods opened the way to sieve theory, leading to upper
bounds of the right order of magnitude for the number of twin primes
$p \leq x$ and to the important achievement of Jingrun Chen \cite{Che} from
1966 that $\#\{p \leq x: p+r =P_2\} \gg \frac{x}{\log^2 x}$, where,
for an integer $k$, $P_k$ denotes the product of at most $k$ primes.
This result relies on another important application of sieve theory,
the Bombieri-Vinogradov theorem on averages of primes in an
arithmetic progression, obtained independently by E. Bombieri and
A.I. Vinogradov  in the mid 1960s. In the late 1980s,  H. Maier and C. Pomerance, and,
subsequently,  the first author of this paper,  obtained similar Bombieri-Vinogradov type results
concerning averages of twin primes (see \cite{MaPo}, \cite{Ba}),  by building on previous work of
N.G. Chudakov,  A.F. Lavrik,  H.L. Montgomery and R.C. Vaughan.

The twin prime conjecture can be generalized in many directions.
For instance, the Hardy-Littlewood heuristics can  be used to
predict the (same) asymptotic formula for the number of primes $p
\leq x$ such that $\frac{p-1}{2}$ is also a prime. This question
may be reformulated as counting the number of primes $p \leq x$
such that the group $\F_p^{\ast} \backslash \{\pm1\}$ is of prime
order. Such a reformulation may then be easily generalized to
other groups, say to the group of points of an elliptic curve:
given an elliptic curve $E/\Q$ over the field of rational
numbers, count the number of primes $p \leq x$ of good reduction
for $E$ such that the group $E(\F_p)/E(\Q)_{\text{tors}}$ is of
prime order, where $E(\F_p)$ denotes the reduction of $E$ modulo
$p$ and $E(\Q)_{\text{tors}}$ denotes the torsion subgroup of
$E/\Q$. This question has theoretical relevance to elliptic curve
cryptography and was first considered by Neal Koblitz in 1988:

\noindent
{\bf{Koblitz's  Conjecture}} \label{koblitz}
 \cite{Ko}

\noindent
{\it{
 Let $E/\Q$ be an elliptic curve
defined over the field of rational numbers. We assume that $E$ is
not $\Q$-isogenous to an elliptic curve with non-trivial
$\Q$-torsion subgroup. Then there exists a positive constant
$C(E)$ such that, as $x \rightarrow \infty$,
$$
\pi_E^{\rm twin}(x) := \# \left\{p \leq x:  |E(\F_p)| \; \text{is
prime}\right\} \sim C(E) \frac{x}{\log^2 x}.
$$
}}

\noindent
A candidate for the explicit constant  $C(E)$ was given by Koblitz
in his paper  and was later corrected by D. Zywina \cite{Zy}. It is
described in detail in Section \ref{constant}  in the generic case
that  $E/\Q$ is without complex multiplication.

It is useful to write the number of points of $E$ over $\F_p$ as
\begin{eqnarray*}
\vert E(\F_p) \vert = p+1-a_p(E),
\end{eqnarray*}
where $a_p(E)$ satisfies the Hasse bound $| a_p(E) | \leq 2
\sqrt{p}$. This also makes the analogy between Koblitz's
Conjecture and the twin prime conjecture more apparent.
Exploiting this analogy, one can employ sieve methods to find
partial results towards Koblitz's Conjecture. This approach was initiated
by S.A. Miri and V.K. Murty \cite{MiMu} and further refined by
A. Steuding and J.  Weng \cite{StWe}, the second
author \cite{Co}, and H. Iwaniec and J. Jimenez Urroz
\cite{IwUr}. We currently know upper bounds of
the right order of magnitude for $\pi_E^{\rm twin}(x)$ (\cite{Co}), provided
the Generalized Riemann Hypothesis  holds if $E/\Q$ is
without complex multiplication,  and various lower bounds in the style of Chen's result (\cite{MiMu},
\cite{StWe}, \cite{Co}, \cite{IwUr}).
 Regarding  lower bounds,
the best result that one may hope to achieve by sieve technology
was obtained by \cite{IwUr} for the complex multiplication elliptic
curve $E: y^2 = x^3 - x$. They showed that
$$
\#
\left\{ p \leq x :  p \equiv 1 (\mod 4), \frac{1}{8} \left\vert E(\F_p)
\right\vert = P_2 \right\} \gg \frac{x}{\log^2{x}}.
$$


The main purpose of this paper is to prove the validity of Koblitz's Conjecture on average over
a set of elliptic curves $E/\Q$:
\begin{theorem}
\label{thm1} Let $x>0$ be a variable and let $A=A(x), B=B(x)$ be
parameters such that $A,B > x^{1/2 +  \varepsilon}$ and $AB >
x^{3/2+\varepsilon}$ for any fixed $\varepsilon > 0$. Let
$\mathcal{C}$ be the set of elliptic curves $E(a, b): \; \; Y^2 = X^3 + aX + b$,
where $a,b \in \bbz$ with $|a| \leq A, |b| \leq B$. Then, as $x
\rightarrow \infty$,
\begin{eqnarray*}
\label{averageTWIN}
\frac{1}{|{\mathcal{C}}|}
\ds\sum_{E \in \mathcal{C}}
\pi_E^{\rm twin}(x)
\sim
\const \frac{x}{\log^2{x}},
\end{eqnarray*}
where $\const$ is the non-zero constant
\begin{eqnarray*}
\label{formulaforC} \const &:=& \frac{2}{3} \ds\prod_{{\ell \neq 2}}
\frac{\ell^4-2\ell^3-\ell^2+3 \ell}{(\ell-1)^3 (\ell+1)} =
\ds\prod_{\ell} \left( 1 - \frac{\ell^2-\ell-1}{(\ell-1)^3(\ell+1)}
\right).
\end{eqnarray*}
\end{theorem}

As will be shown in Section \ref{constant}, the average constant
$\mathfrak{C}$ gives further evidence for the conjectural constant
of Koblitz's Conjecture.

The first steps in the proof of Theorem \ref{thm1} follow the ones in the proof of the average
Lang-Trotter Conjecture on Frobenius traces obtained  by the third author and F. Pappalardi
in  \cite{DaPa},  with an improvement due to  S. Baier \cite{Bai} to shorten the average in terms
of $A$ and $B$ (in \cite{DaPa}, one needed to take $A, B > x^{1+\varepsilon}$).
More precisely, we first reduce the average of Koblitz's Conjecture given in  Theorem \ref{thm1}
to an average involving
only elliptic curves over the finite field $\F_p$; we then use Deuring's Theorem to rewrite this as an
average of Kronecker class numbers.
These steps are described in detail in Section \ref{reducetoH}.

As corollaries,  we obtain  the average result of
Theorem \ref{thm1}, as well as:

\begin{theorem} \label{thm2}
Let
$$
\pi^*(p) := \# \left\{E/\F_p:  |E(\F_p)| \;\mbox{is
prime} \right\}.
$$
Then, as $x \rightarrow \infty$,
\begin{eqnarray*}
\ds\sum_{p \leq x} \pi^*(p) \sim \frac{ \const x^3}{3 \log^2{x}},
\end{eqnarray*}
where $\mathfrak{C}$ is the constant of Theorem \ref{thm1}.
\end{theorem}

One remarks that
Theorem \ref{thm2} is only concerned with  the
distribution of elliptic curves over $\F_p$ having  a prime number of
points. More properties of this distribution could be obtained by
considering the higher moments
$$
M_k(x) := \ds\sum_{p \leq x} \left( \pi^*(p) \right)^k \; \; \text{for $k \geq 1$}.
$$

One of the key ingredients in the proof of Theorems \ref{thm1}-\ref{thm2}
is an average of the standard twin prime conjecture. Such
averages were first considered by N.G. Chudakov  \cite{Chu} and A.F.  Lavrik \cite{Lav},  and then,
among others, by A. Balog \cite{Ba}, who added distribution in residue classes, and by
A. Perelli and J. Pintz \cite{PePi}, who shortened the average. The
length of the average needed for our application is dictated by
Hasse's bound and is {\it{short}} ($\sqrt{x}$ compared to $x$); additionally,  we
also need distribution in {\it{residue classes}}. Such a  mixture
of additional features  is not in the literature and is proven in our paper:

\begin{theorem} \label{thmAB}
Let $x > 0$ and let $\varepsilon, M > 0$. Then  there exists an integer
$N(M) > 0$ such that,  for any $x^{1/3+\varepsilon} \leq R
\leq x$, $N > N(M)$, $Q \leq x {\log^{-N}{x}}$,
and $X, Y$ satisfying $X+Y \leq x$, we have
$$
\ds\sum_{0<|r| \leq R} \ds\sum_{q \leq Q} \ds\sum_{a (\mod q)}
\left\vert \ds\sum_{{{X < p \leq X+Y} \atop {p \equiv a (\mod
q)}} \atop{p-p'=r}} \log{p} \cdot \log{p'} - {\mathfrak S}(r,q,a)
Y \right\vert^2 \ll \frac{R x^2}{\log^{M}{x}},
$$
where
\begin{eqnarray*}
{\mathfrak S}(r,q,a)
&:=&
\left\{ \begin{array}{cc}
\displaystyle{\frac{1}{\phi(q)} {\mathfrak S}(rq)} & \mbox{if $2
\mid r, (a,q)=(a-r,q)=1$},
\\ 0 & \mbox{otherwise},
\end{array} \right.
\end{eqnarray*}
and $\mathfrak S(r q)$ is as in (\ref{twin-prime-constant}) and $\phi(q)$ is the Euler function of
$q$.
Here (and in what follows), $q$ denotes positive integers, and $p, p'$ and $\ell$ denote
rational primes.
\end{theorem}

The structure of the  paper is as follows. In Section \ref{constant}, we present the heuristic
reasoning behind Koblitz's Conjecture and discuss the constant $\mathfrak{C}$.
In Section \ref{reducetoH}, we reduce the statements of Theorems \ref{thm1}-\ref{thm2} to an average of
Kronecker class numbers (Proposition \ref{thmwithXY}).
In Section \ref{reducetoTWIN},  we show how an average of the twin prime conjecture implies
Proposition \ref{thmwithXY}. Finally, in Sections \ref{ATPC}-\ref{main}, we give the proofs of the
afore-mentioned average of the twin prime conjecture and of Theorem \ref{thmAB}.

\section{Average of Koblitz's Conjecture and the conjectural constant}
\label{constant}

\hspace*{0.5cm}
The constant $C(E)$ of  Koblitz's Conjecture  is based on the
following heuristic argument, which is reminiscent of the argument
leading to the classical twin prime constant of Hardy and Littlewood
(see, for example, \cite{So}).

We want to count the number of primes
$p$ such that $p+1-a_p(E)$ is also a  prime. For each prime $\ell$, this
means that $p+1-a_p(E)$ is not divisible by $\ell$. For a random
integer $n$, the probability that $\ell \nmid n$ would be
$(\ell-1)/\ell$. To compute the probability that $\ell \nmid
p+1-a_p(E)$ we consider the action of $\mbox{Gal}(\bar{\Q}/\Q)$ on
the group  $E[\ell]$ of $\ell$-torsion points of $E$, which leads to the
injection
$$
\rho_\ell: \mbox{Gal}(\Q(E[\ell])/ \Q) \hookrightarrow  \mbox{GL}_2
(\Z/\ell\Z);
$$
here,  $\Q(E[\ell])$ is the field obtained by adjoining to $\Q$  the
coordinates of the points in $E[\ell]$.
By studying the action of the Frobenius map on the torsion points of
$E$, it follows  that
\begin{eqnarray*}
\mbox{tr}(\rho_\ell(\sigma_p)) &\equiv& a_p(E) (\mod \ell), \\
\det(\rho_\ell(\sigma_p)) &\equiv& p (\mod \ell),
\end{eqnarray*}
for all primes $p \neq \ell$ of good reduction for $E$.
Then  the probability that $\ell \nmid p+1-a_p(E)$ can  be
evaluated by counting matrices $g$ in $\mbox{GL}_2(\Z / \ell \Z)$
such that $\det(g) + 1 - \mbox{tr}(g) \not\equiv 0 (\mod \ell)$.

Let $G(\ell)$ be the image of $\rho_\ell$ in $\mbox{GL}_2
(\Z/\ell\Z)$ and let
\begin{eqnarray*}
\Omega(\ell) &:=& \left\{ g \in G(\ell):  \det(g)+1-\tr(g) \equiv 0 (\mod \ell) \right\}, \\
\Omega^\prime(\ell) &:=& \left\{ g \in \mbox{GL}_2(\Z/\ell\Z):
\det(g)+1-\tr(g) \equiv 0 (\mod \ell) \right\}.
\end{eqnarray*}
Then, at each prime $\ell$, the correcting probability factor is the
quotient
$$
\frac{\displaystyle{1 -
\frac{|\Omega(\ell)|}{|G(\ell)|}}}{\displaystyle {1 -
\frac{1}{\ell}}},
$$
where the numerator is the probability that $p+1-a_p(E)$ is not
divisible by $\ell$ and the denominator is the probability that a
random integer is not divisible by $\ell$.

If $G(\ell)=\mbox{GL}_2(\Z/\ell\Z)$, then we have
\begin{eqnarray*}
\frac{\displaystyle{1 -
\frac{|\Omega(\ell)|}{|G(\ell)|}}}{\displaystyle {1 -
\frac{1}{\ell}}} = \frac{\displaystyle{1 -
\frac{|\Omega^\prime(\ell)|}{|\mbox{GL}_2(\Z/\ell\Z)|}}}{\displaystyle
{1 - \frac{1}{\ell}}}
 &=&
1 - \frac{\ell^2-\ell-1}{(\ell-1)^3 (\ell+1)}.
\end{eqnarray*}

The constant $C(E)$ of \cite{Ko} is defined as the product over all
primes $\ell$ of the local factors above. In \cite{Zy}, Zywina made
the observation that the probabilities are not independent from one
prime to another, because the fields $\Q(E[\ell])$
are never independent for all primes $\ell$,  as already
observed by Serre in \cite{Se}. Their  dependence can be quantified:
for each elliptic curve $E/\Q$,  there is an integer $M_E$ which has the
property that the probabilities are independent for primes $\ell
\nmid M_E$ (see \cite{Zy} for a precise definition).

Now  let
$G(M_E) \subseteq \GL_2(\Z/M_E\Z)$ be the Galois group of
$\Q(E[M_E])/\Q$ and let
$$
\Omega(M_E) := \left\{ g \in G(M_E):  \det(g)+1-\tr(g) \equiv 0 (\mod {M_E}) \right\} .
$$
For each elliptic curve $E$, the constant $C(E)$ of  Koblitz's Conjecture is then
expected  to be
\begin{eqnarray} \label{explicitKC}
C(E) &=& \frac{1 -
\displaystyle{\frac{|\Omega(M_E)|}{|G(M_E)|}}}{\;\displaystyle{\prod_{\ell
| M_E}} \displaystyle{1-\frac{1}{\ell}}}
 \times \; \prod_{\ell \nmid M_E} \left(1 - \frac{\ell^2-\ell-1}{(\ell-1)^3
 (\ell+1)}\right).
\end{eqnarray}

We remark that even if $M_E$ is never 1, it is possible for some
elliptic curves  $E/\Q$ to have
$$
C(E) = \prod_{\ell} \left(1 - \frac{\ell^2-\ell-1}{(\ell-1)^3
 (\ell+1)}\right) = \mathfrak{C},
$$
as  shown in \cite{Jo} and \cite{Zy} (it is shown in \cite{Jo}
that $C(E) = \mathfrak{C}$ when the square-free part of the discriminant 
of $E$ is congruent to 2 or 3 modulo 4).  However,  the average constant $\mathfrak{C}$
of Theorem \ref{thm1} should not be thought of as the constant of {\it{any}}
given curve over $\Q$, but as the {\it{average}} of all the constants
$C(E)$. Indeed, in \cite{Jo}, N.  Jones shows that
if one assumes  a positive answer to a well-known question of Serre regarding the
open image theorem for elliptic curves proven in \cite{Se},
then the average of the conjectural constants
$C(E)$ of (\ref{formulaforC}) is indeed the average constant of
Theorem \ref{thm1}. Our result then gives evidence for both the
asymptotic of Koblitz's Conjecture and  the constant
appearing in the conjecture.


\section{Reduction to an average of Kronecker class numbers and proofs of
Theorems \ref{thm1}-\ref{thm2}}

\label{reducetoH}

 \hspace*{0.5cm}
 Let $x>0$ and let $\mathcal{C}$ be
the family of elliptic curves introduced in Theorem \ref{thm1}. In
this section we show how the average of $\pi^{\rm{twin}}_E(x)$ over
$E/\Q$ reduces to an average of  Kronecker class numbers, which is
Proposition \ref{thmwithXY} of the present section.

First we write
\begin{eqnarray} \label{firstreduction}
\frac{1}{\vert \mathcal{C} \vert} \ds\sum_{E \in \mathcal{C}}
\pi^{\rm{twin}}_E(x) &=& \frac{1}{\vert \mathcal{C} \vert}
\ds\sum_{{p \leq x, \; |r| \leq 2 \sqrt{p}} \atop {p+1-r
\;\rm{prime}}} N_{A,B,r}(p),
\end{eqnarray}
where
$$ N_{A,B,r}(p) :=  \# \left\{ |a| \leq A, |b| \leq B: a_p(E(a,b)) = r \right\}.$$
As the value of $a_p(E(a,b))$ depends
only on $a$ and $b$ modulo $p$, $N_{A,B,r}(p)$ can be evaluated by
counting elliptic curves over the finite field $\F_p$ and having  $p+1-r$
points for each given $r$. In turn, this can be done  using the following results:

\begin{theorem}(Deuring's Theorem \cite{De})

\noindent Let $p > 3$ be a prime and let $r$ be an integer such that
$r^2 - 4p < 0$. Let $\mathcal{E}_r(p)$ be the set of
$\F_p$-isomorphism classes of elliptic curves over $\F_p$ having
$p+1-r$   $\F_p$-rational points. Then
$$
\ds\sum_{E \in \mathcal{E}_r(p)} \frac{1}{\# \Aut(E)} = H(r^2 - 4p),
$$
where $\Aut(E)$ is the automorphism group of $E$ and for any $D<0$,
$H(D)$ is the Kronecker class number
$$
H(D) := \sum_{{f^2 \mid D} \atop {\frac{D}{f^2} \equiv 0,1 (\mod 4)}}
\frac{h(D/f^2)}{w(D/f^2)}
$$
defined in terms of the class number $h(D/f^2)$ and number of
units $w(D/f^2)$ of $\Q(\sqrt{D/f^2})$.\\
Then,  for any fixed $-2 \sqrt{p} \leq r \leq 2 \sqrt{p}$, there
are exactly $(p-1) H(r^2-4p)$ elliptic curves defined over $\F_p$
with $p+1-r$ points.
\end{theorem}

\begin{lemma}\label{trivialbound} Let $D$ be a positive integer such that
$-D \equiv 0,1 (\mod{4})$. Then, as $D \rightarrow \infty$,
\begin{eqnarray}
\label{boundH} H(-D) \ll D^{1/2} \log^2{D}.
\end{eqnarray}
\end{lemma}
\begin{proof} This follows from the class number
formula and from standard bounds on special values of Dirichlet
L-functions.\end{proof}

By reducing modulo $p$ and counting curves
over $\F_p$,  it follows that
\begin{eqnarray} \label{likeDaPa}
N_{A,B,r}(p) &=& \left( \frac{2A}{p} + \O(1) \right) \left(
\frac{2B}{p} + \O(1) \right) \left( p H(r^2-4p) + \O\left(\sqrt{p} \log^2{p}\right) \right).
\end{eqnarray}

This is the approach of \cite{DaPa}  which  gives an average of size
$A, B > x^{1+\varepsilon}$ when one replaces (\ref{likeDaPa}) in
(\ref{firstreduction}). In order to get the shorter average of
Theorem \ref{thm1}, one needs to count $N_{A,B,r}(p)$ in a more subtle
way, using character sums. This was done by Baier,  who
proved the following lemma:

\begin{lemma} \label{baier} Let $N_{A,B,r}(p)$ be the number of
curves $E(a,b) \in \mathcal{C}$ such that $a_p(E(a,b))=r$. Then,
\begin{eqnarray*}
N_{A,B,r}(p) &=& \frac{4AB \cdot H(r^2-4p)}{p} + \O \left( \frac{AB}{p} +
\frac{AB \cdot H(r^2-4p)}{p^2} + A + B \right. \\
&&
\left. + (AB \cdot H(r^2-4p))^{1/2} \log^3{p} +
\frac{(A+B) H(r^2-4p)}{p^{1/2}} \log{p} \right).
\end{eqnarray*}
\end{lemma}
\begin{proof}
See \cite{Bai}.
\end{proof}

Replacing Lemma \ref{baier} in (\ref{firstreduction}) and using
Lemma \ref{trivialbound} to bound the Kronecker class numbers, we
obtain that
\begin{eqnarray} \nonumber \frac{1}{\vert \mathcal{C} \vert} \ds\sum_{E \in \mathcal{C}}
\pi^{\rm{twin}}_E(x) &=& \sum_{{p \leq x, \; |r| \leq 2 \sqrt{p}}
\atop {p+1-r \;\rm{prime}}} \frac{H(r^2-4p)}{p} + \O \left( x^{1/2} +
\frac{x^{3/2} \log^3{x}}{A} + \frac{x^{3/2} \log^3{x}}{A} +
\frac{x^{7/4} \log^4{x}}{(AB)^{1/2}} \right) \\
\label{forthm1} &=& \sum_{{p \leq x, \; |r| \leq 2 \sqrt{p}} \atop
{p+1-r \;\rm{prime}}} \frac{H(r^2-4p)}{p} + \O(x^{1-\varepsilon}),
\end{eqnarray}
provided we take $A, B$ such that $A,B > x^{1/2 + 2 \varepsilon}$ and $AB >
x^{3/2+2\varepsilon}$.

In a similar way,  using Deuring's Theorem and the bound of Lemma
\ref{trivialbound}, we can write
\begin{eqnarray} \label{forthm2}
\sum_{p \leq x} \pi^*(p) &=&  \sum_{{p \leq x, \; |r| \leq 2
\sqrt{p}} \atop {p+1-r \;\rm{prime}}} p H(r^2-4p) + \O\left( x^2
\log^2 x\right).
\end{eqnarray}
Thus, once again,  we need to evaluate an average of class numbers.
Then, Theorems \ref{thm1} and \ref{thm2} will follow from:

\begin{proposition}
\label{thmwithXY} Let $x, X, Y$ be positive real numbers such that
$X+Y \leq x$, $Y \geq \sqrt{X}$.
Then, for any $M > 0$,
\begin{eqnarray*}
\ds\sum_{X < p \leq X+Y} \ds\sum_{{|r| \leq 2 \sqrt{X}} \atop
{p+1-r \;\rm{prime}}} p H(r^2-4p) &=& \frac{\const
X^2Y}{\log^2{X}} + \O \left( XY^2 \log{X} \right) + \O \left(
\frac{x^3}{\log^M{x}} \right)
\end{eqnarray*}
where $\const$ is the constant of Theorem \ref{thm1}.
\end{proposition}

Now  let us indicate in detail  how Proposition \ref{thmwithXY} implies Theorems
\ref{thm1} and  \ref{thm2}.

\noindent
{\bf{Proof of Theorems \ref{thm1}-\ref{thm2}.}}
Comparing (\ref{forthm1}) and
(\ref{forthm2}), we see that the two asymptotics to prove are
equivalent by partial summation; we will only prove the second one.

Let $M$ be any positive integer (e.g.  $M=10$ suffices) and
$K := [\log^{M/2}x]$.
Let
$$
Y := \frac{x}{K} \; \;  \text{and} \; \; X  := k Y  \; \; \text{for} \;  \; 0 \leq k \leq K-1.
$$
We partition the interval $p \leq x$ into $K$
intervals of length $Y$ and rewrite the main term of
(\ref{forthm2}) as
\begin{eqnarray}
\nonumber \ds\sum_{ {{p \leq x} \atop {|r| \leq 2 \sqrt{p}}}
\atop {p+1-r \;\rm{prime}} } p H(r^2 - 4p)
&=&
\ds\sum_{0 \leq k \leq K-1}
\ds\sum_{ {{X < p \leq X+Y} \atop {|r| \leq 2
\sqrt{p}}} \atop {p+1-r \;\rm{prime}} } p H(r^2-4p)
\\
&=&
\ds\sum_{0 \leq k \leq K-1}
\ds\sum_{ {{X < p \leq X+Y}
\atop {|r| \leq 2 \sqrt{X}}} \atop {p+1-r \;\rm{prime}} } p H(r^2-4p) \label{withET}
+
\O \left( \ds\sum_{0 \leq k \leq K-1} X Y^2 \log^2{x} \right),
\end{eqnarray}
where the $\O$-term comes from the bound (\ref{boundH}) on the
Kronecker  class number $H(r^2-4p)$ for $r$ in the interval $2
\sqrt{X} < r \leq 2 \sqrt{p}$ and is bounded by $K^2Y^3\log^2x\leq
x\log^{2-M/2}x$. For the main term we use Propostion \ref{thmwithXY}
with the same $M$, and obtain:
\begin{eqnarray}
\nonumber
\ds\sum_{0 \leq k \leq K-1}
\ds\sum_{ {{X < p \leq X+Y} \atop {p+1-r \;\rm{prime}}} \atop {|r| \leq 2\sqrt{X}}}
p H(r^2-4p)
&=&
\const Y^3 \ds\sum_{1 \leq k \leq K-1}
\frac{k^2}{\log^2(kY)} + \O\left(\ds\sum_{1 \leq k \leq K-1} X Y^2 \log X\right)
\\
\nonumber
&& +
\O\left(\ds\sum_{1 \leq k \leq K-1} \frac{x^3}{\log^{M} x}\right) + \O\left(Y^{5/2} \log^2 Y\right)
\\
\nonumber
&=&
\const Y^3 \int_{1}^{K-1} \frac{t^2}{\log^2(tY)} \; d t + \O\left( \frac{x^3}{\log^{M/2-1}{x}} \right)
\\
\nonumber &=& \const \int_{Y}^{x} \frac{u^2}{\log^2(u)} \; d u +
\O\left( \frac{x^3}{\log^{M/2-1}{x}} \right)
\\
&=& \frac{\const x^3}{3 \log^2{x}} + \O \left(
\frac{x^3}{\log^3{x}} \right).\label{asymptotic}
\end{eqnarray}
Replacing(\ref{withET}) and (\ref{asymptotic}) in (\ref{forthm2}),
the proof of Theorem \ref{thm2} is completed. The proof of Theorem
\ref{thm1} follows by partial summation.
\qed


\section{Reduction to an average twin prime conjecture and proof of Proposition \ref{thmwithXY}}

\label{reducetoTWIN}

\hspace*{0.5cm}
In this section we show how Proposition \ref{thmwithXY} reduces to an average of the
twin prime conjecture which, in turn,  will be proved completely in Sections \ref{ATPC}-\ref{main}.
To be precise,
our proof of Proposition \ref{thmwithXY} relies on the validity of the following result:
\begin{proposition}
\label{mainone}
Let $x > 0$ and let $M > 0$. Let $X, Y, R,
U, V$ be parameters depending on $x$ and satisfying
$$
X+Y \leq x, \; R \leq x, \; x^{1/2} \log^N{x} \leq U, \;
\log^N x \leq V, \; UV^2\leq x\log^{-N}x.
$$
There exists an integer $N(M) > 0$ such that,
if  $N>N(M)$, then, as $x \rightarrow \infty$,
\label{wantedATP}
\begin{eqnarray}
&&\nonumber \ds\sum_{{{|r| \leq R} \atop {f \leq V}} \atop {n \leq
U}} \frac{1}{nf} \ds\sum_{a (\mod 4n)} \left( \frac{a}{n} \right)
\ds\sum_{{{X < p \leq X+Y}\atop p+1-r \text{ prime}}\atop p\equiv
(r^2-af^2)/4\, (\mod nf^2) }\log{p} \cdot \log{(p+1-r)}
\\ \label{EQwantedATP} &&= 2 \const RY + \O \left(
\frac{Rx}{\log^M{x}}+x^{4/3+\varepsilon} \right),
\end{eqnarray}
where $\const$ is the  constant of Theorem \ref{thm1}.
\end{proposition}
We assume this result as true and proceed to proving Proposition \ref{thmwithXY}.

\noindent {\bf{Proof of Proposition \ref{thmwithXY}}.}
Let $x, X, Y$
be as in the statement of Proposition \ref{thmwithXY}. Using the
class number formula, we write
\begin{eqnarray}
\label{eqbefore} \ds\sum_{X < p \leq X+Y} \ds\sum_{{|r| \leq 2
\sqrt{X}} \atop {p+1-r \;\rm{prime}}} p H(r^2-4p) &=& \frac{1}{2
\pi} \ds\sum_{{|r| \leq 2\sqrt{X}} \atop {f \leq 2\sqrt{x}}}
\frac{1}{f} {\ds\sum_{X < p \leq X+Y}}^{*} p {\sqrt{4p-r^2}} L(1,
\chi_d),
\end{eqnarray}
where $d = d(r,p,f) := (r^2 - 4p)/f^2$ and the $*$ on the summation over $p$
indicates that we are summing over primes
$X < p \leq X+Y$ such that
$$
p+1-r \;{\rm{prime}},
\;f^2 \mid r^2-4p,
\; \text{and} \; d \equiv 0,1 (\mod 4).
$$
Here, $\chi_d$ denotes the Dirichlet character modulo $d$ defined
by the Kronecker symbol and $L(s, \chi_d)$ denotes its Dirichlet
L-function.

Using the P\'olya-Vinogradov inequality, we write the special
value $L(1, \chi_d)$ as
\begin{eqnarray*}
L(1, \chi_d) &=& \sum_{n \leq U} \frac{\chi_d(n)}{n} + \sum_{n > U} \frac{\chi_d(n)}{n}\\
&=& \sum_{n \leq U} \frac{\chi_d(n)}{n} + \O \left( \frac{\sqrt{|d|}\log{|d|}}{U} \right),
\end{eqnarray*}
where $U=U(x)$ is a parameter to be chosen soon.
By plugging  this into (\ref{eqbefore}), we obtain that
\begin{eqnarray}
\nonumber && \ds\sum_{X < p \leq X+Y} \ds\sum_{{|r| \leq 2
\sqrt{X}} \atop {p+1-r \;\rm{prime}}} p H(r^2-4p)
\\
\label{last} && \hspace{1cm} = \frac{1}{2 \pi} \ds\sum_{ {{|r|
\leq 2\sqrt{X}} \atop {f \leq 2\sqrt{x}} } \atop {n \leq U}}
\frac{1}{nf} {\ds\sum_{X < p \leq X+Y}}^* p {\sqrt{4p-r^2}}
\chi_d(n) + \O \left( \frac{x^{7/2}}{U} \right).
\end{eqnarray}
Thus, by taking
\begin{eqnarray}
\label{conditiononU-1}
 x^{1/2} \log^{M + 1}{x} \leq U \leq x,
\end{eqnarray}
the $\O$-term above becomes
$\O \left( x^3/\log^M x \right)$.

Now let us also truncate the sum over $f$. We write
\begin{eqnarray}
\nonumber && \frac{1}{2 \pi} \ds\sum_{ {{|r| \leq 2\sqrt{X}}
\atop {f \leq 2\sqrt{x}} } \atop {n \leq U}} \frac{1}{nf}
{\ds\sum_{X < p \leq X+Y}}^* p {\sqrt{4p-r^2}} \chi_d(n)
\\
\nonumber && \hspace{1cm} = \frac{1}{2 \pi} \ds\sum_{ {{|r| \leq 2\sqrt{X}}
\atop {f \leq V} } \atop {n \leq U}} \frac{1}{nf} {\ds\sum_{X < p
\leq X+Y}}^* p {\sqrt{4p-r^2}} \chi_d(n) + \O\left( \ds\sum_{
{{|r| \leq 2\sqrt{X}} \atop {V < f \leq 2\sqrt{x}} } \atop {n
\leq U}} \frac{1}{nf} \ds\sum_{ X < p \leq X + Y \atop{p \equiv
\frac{r^2}{4} (\mod f^2)} } p^{3/2} \right)
\\
\label{last-V} && \hspace{1cm} = \frac{1}{2 \pi} \ds\sum_{ {{|r| \leq
2\sqrt{X}} \atop {f \leq V} } \atop {n \leq U}} \frac{1}{nf}
{\ds\sum_{X < p \leq X+Y}}^* p {\sqrt{4p-r^2}} \chi_d(n) +
\O\left(\frac{ x^{3} \log U}{V^2}\right).
\end{eqnarray}
On the second line we  used that since $r$ is odd and $f^2
| r^2 - 4p$, we must have that $f$ is odd, hence the condition in
the sum over $p$ that $4 p \equiv r^2 (\mod f^2)$ becomes $p \equiv
\bar{4} {r^2} (\mod f^2)$, where $\bar{4}$ is the inverse of 4
modulo $f^2$. The other estimates used for the $\O$-term are
elementary.

We choose $V$ such that
\begin{eqnarray}
\label{condition-on-V}
V \geq \left(\log{x} \right)^{(M+1)/2}
\end{eqnarray}
and then the $\O$-term above becomes $\O \left( x^3 / \log^M{x} \right)$.

Now we use quadratic reciprocity and
consider $\chi_d(n)$ as a character modulo  $4n$. Hence we rewrite the main term of (\ref{last-V}) as
$$
\frac{1}{2\pi} \ds\sum_{{{|r| \leq 2\sqrt{X}} \atop {f \leq V}}
\atop {n \leq U}} \frac{1}{nf} \ds\sum_{a (\mod 4n)} \left(
\frac{a}{n} \right) {\ds\sum_{X < p \leq X+Y}}^{**} p
\sqrt{4p-r^2},
$$
where the $**$ on the summation over $p$
indicates that we are  summing over primes $X < p \leq X+Y$ such that
$$
 p+1-r \;{\rm{prime}},
\; f^2 \mid r^2-4p, \; d \equiv 0,1 (\mod 4), \; \text{and} \;
\frac{r^2 - 4p}{f^2} \equiv a (\mod 4n).
$$
Since $r$ and $f$ must be odd, we necessarily have   $d \equiv  1
(\mod 4)$; thus $**$ is equivalent to the conditions
$$
p+1-r \;{\rm{prime}}
\; \text{and} \;
 p \equiv \frac{r^2 - af^2}{4} (\mod nf^2).
$$

We now change the weight of the primes $p$ from $p \sqrt{4p-r^2}$ to
$$
\frac{\log{p} \cdot \log{(p+1-r)}}{\log^2{X}} X \sqrt{4X-r^2}.
$$
Then the main term of (\ref{last-V}) becomes
\begin{eqnarray}
\label{tag1} && \frac{X}{2 \pi \log^2{X}} \ds\sum_{ {{|r| \leq
2\sqrt{X}} \atop {f \leq V}} \atop {n \leq U}} \frac{1}{nf}
\sqrt{4X-r^2} \ds\sum_{a (\mod 4n)} \left( \frac{a}{n} \right)
{\ds\sum_{X < p \leq X+Y}}^{**} \log{p} \cdot \log{(p+1-r)}
\\
\nonumber
&&
\hspace{1in} + \O \left( X Y^2 \log{X} \right)
\end{eqnarray}
and so we reduced our question to an average of the standard twin prime conjecture,
twisted by some Kronecker symbols. This average is evaluated using Proposition \ref{mainone}
stated in the beginning of this section; the details follow.


Let us write the left hand side of Proposition
\ref{mainone} as $\ds\sum_{|r| \leq R} F(r)$. With this notation,
the main term of (\ref{tag1}) becomes
\begin{eqnarray*}
&& \frac{X}{2 \pi \log^2{X}} \ds\sum_{|r| \leq 2\sqrt{X}} F(r)
\sqrt{4X-r^2},
\end{eqnarray*}
which  we can compute from Proposition \ref{mainone} by partial
summation. We obtain:
\begin{eqnarray*}
&& \frac{X}{2 \pi \log^2{X}} \ds\sum_{|r| \leq 2\sqrt{X}} F(r)
\sqrt{4X-r^2}
\\
&&= \label{tag2} \frac{X}{2 \pi \log^2{X}} \int_0^{2\sqrt{X}}
\left( 2 \const t Y + \O \left(
\frac{tx}{\log^{M}{x}}+x^{4/3+\varepsilon} \right) \right)
\left( t (4X-t^2)^{-1/2} \right) \;dt
\\
&&= \frac{\const}{\pi} \cdot \frac{XY}{\log^2 X} \int_0^{2
\sqrt{X}} t^2 (4X -t^2)^{-1/2} \;dt + \O\left(\frac{(Xx)^{3/2}}{\log^2
X\log^M x}\right).
\end{eqnarray*}
Evaluating the integral
$$
\int_0^{2\sqrt X} t^2 (4X-t^2)^{-1/2} \;dt = 4X\int_0^1
t^2(1-t^2)^{-1/2}\,dt = \pi X
$$
and replacing in the above, we obtain
\begin{eqnarray*}
\frac{X}{2 \pi \log^2{X}} \ds\sum_{|r| \leq 2\sqrt{X}} F(r)
\sqrt{4X-r^2} &=& \frac{\const X^2 Y}{\log^2{X}} +
\O\left(\frac{(X x)^{3/2}}{(\log^2 X) (\log^M x)}\right).
\end{eqnarray*}
Using this together with (\ref{tag1}), the proof of Proposition
\ref{thmwithXY}  is now completed
(provided that Proposition \ref{mainone} holds).
\qed


\section{Average of the twin prime conjecture and proof of Theorem \ref{thmAB}}
\label{ATPC}

\hspace*{0.5cm}
In this section  we shall  prove Theorem \ref{thmAB}. The statement is a
Barban-Davenport-Halberstam type distribution result for twin
primes, where the average is  over the twin prime differences. The main
(and difficult)  part of the proof is
the  case $Q=1$ of Theorem \ref{thmAB}.  A version of this was proven by Perelli and Pintz in
\cite{PePi}. Beside minor cosmetics, their result differs from
what we need in two aspects: rather than a Goldbach type problem,
we have a twin prime problem; more importantly,  we need a
Siegel-Walfisz type analogue, namely:

\begin{proposition}
\label{SiegelWalfisz}
Let $\varepsilon, M, N >0$ be fixed.
Then there exists $x(\varepsilon,M,N) > 0$ such that,
for any $x>x(\varepsilon,M,N)$, $x^{1/3+\varepsilon}\leq R\leq x$,
$q\leq\log^Nx$, $(a,q)=1$ and $0<X<X+Y\leq x$,  we have
$$
\sum_{0<r\leq R}\left| \sum_{X<p\leq X+Y \atop{p\equiv a(\mod q)
\atop{p-p'=r}}} \log p \cdot \log p'-\mathfrak{S}(r,q,a)Y\right|^2\ll\frac{Rx^2}{\log^Mx}.
$$
\end{proposition}

The proof of Proposition \ref{SiegelWalfisz} is similar to
the proof of  \cite[Theorem 1]{PePi}. For this reason,  we will only indicate the major
steps that should enable  the interested reader  to modify \cite{PePi}
accordingly. After we complete the proof of Proposition \ref{SiegelWalfisz},
we proceed to  proving  Theorem \ref{thmAB}.

\bigskip\noindent {\bf{Proof of Proposition \ref{SiegelWalfisz}.}}
We use the following notation (which is, unfortunately,  not exactly the
same as in \cite{PePi}):

$$
S_1(\alpha) :=  \ds\sum_{X<p\leq X+Y \atop{p\equiv a(\mod q)}}
\log p \cdot e(p\alpha),
\quad
S_2(\alpha) := \sum_{p'\leq x} \log p' \cdot e(p'\alpha),
\quad
e(y) := e^{2\pi iy},
$$
$$
C := C(\varepsilon, M,N),
\quad
I_{s,b} :=
\text{Farey arc around }\frac bs
= \left\{\frac bs+\eta,\,|\eta|<\frac{\log^{2C}x}{sx}\right\},
$$
$$
\mathfrak{M} := \bigcup_{s\leq\log^Cx}\bigcup_{(b,s)=1}I_{s,b}\,,
\quad
\mathfrak{m} := [0,1]\setminus\mathfrak{M}.
$$

By the circle method, we have
$$
\ds\sum_{X<p\leq X+Y \atop{p\equiv a(\mod q) \atop{p-p'=r}}}
 \log p \cdot \log p'
 =
 \int_0^1 S_1(\alpha)S_2(-\alpha)e(-r\alpha)\,d\alpha,
$$
and so
\begin{eqnarray}
\nonumber
&&\ds\sum_{0<r\leq R}
\left|
\ds\sum_{X<p\leq X+Y \atop{p\equiv a(\mod q) \atop{p-p'=r}}}
\log p\cdot \log p'-\mathfrak{S}(r,q,a)Y
\right|^2 \ll
\ds\sum_{0<r\leq R}\left|\int_{\mathfrak{m}}
S_1(\alpha)S_2(-\alpha)e(-r\alpha)\,d\alpha\right|^2\\
&& \hspace{1in} +
\label{my25}
\ds\sum_{0<r\leq R}\left|\int_{\mathfrak{M}}
S_1(\alpha)S_2(-\alpha)e(-r\alpha)\,d\alpha-\mathfrak{S}
(r,q,a)Y\right|^2.
\end{eqnarray}

The estimate for the contribution of the minor arcs -- the first
term in formula (\ref{my25}) -- is identical to the one
in \cite{PePi}.
To start, we remark that
the Cauchy-Schwarz inequality and the well-known
estimate $\sum_{0<r\leq R} e(ry)\ll \min(R,1/\|y\|)$ reduce this
term to
\begin{eqnarray*}
&&\ds\sum_{0<r\leq R}
\left|
\int_{\mathfrak m} S_1(\alpha)S_2(-\alpha)e(-r\alpha)\,d\alpha
\right|^2
\\&& \hspace{1cm}=
\ds\sum_{0<r\leq R}
\int_{\mathfrak m} S_1(\alpha)S_2(-\alpha)e(-r\alpha)\,d\alpha
\int_{\mathfrak m} \overline{S_1(\beta)S_2(-\beta)}e(r\beta)\,d\beta\\ && \hspace{1cm}\ll
\int_{\mathfrak m}|S_1(\beta)S_2(\beta)|
\int_{\mathfrak m}|S_1(\alpha)S_2(\alpha)|
\min \left(R,\frac1{\|\alpha-\beta\|}\right)\,d\alpha\,d\beta\\ && \hspace{1cm} \ll
\sup_{\beta\in\mathfrak m}\left(\int_{\mathfrak m}
|S_2(\alpha)|^2\min\left(R,\frac1{\|\alpha-\beta\|}\right)^2d\alpha\right)^{1/2}\\
&& \hspace{1cm} \times \left(\int_{\mathfrak m}
|S_1(\alpha)|^2d\alpha\right)^{1/2}\left(\int_{\mathfrak m}
|S_1(\beta)|^2d\beta\right)^{1/2}\left(\int_{\mathfrak m}
|S_2(\beta)|^2d\beta\right)^{1/2},
\end{eqnarray*}
where $\| \|$ denotes the distance to the nearest integer.
Now let us observe that  our $S_2(\alpha)$ is exactly the same as the one  in \cite{PePi},
thus the third integral above can be estimated as in \cite[Section 5]{PePi}.
Note that the function $S_2(\alpha)$  plays the crucial role, while the somewhat different
$S_1(\alpha)$ only appears in Parseval's identity. Since our $S_1(\alpha)$ has
smaller $L^2$-norm than its analogue in \cite{PePi},  the arguments in\cite[Section 3]{PePi} provide
the necessary bound in our case as well.

For the calculation of the major arcs -- the second term in
(\ref{my25}) --  we follow the exact steps of  \cite[Section 4]{PePi}.
First, for $\alpha = \frac{b}{s} + \eta \in I_{s,b}$,
we use the Siegel-Walfisz theorem   to approximate the function
$$
S_1(\alpha)
=
\sum_{X<p\leq X+Y \atop{p\equiv a(\mod q)}} \log p\cdot
e\left(p\left(\frac bs+\eta\right)\right)
=
\sum_{1 \leq c \leq s} e\left(\frac{bc}s\right)
\sum_{X<p\leq X+Y\atop{p\equiv a(\mod q)\atop{p\equiv c(\mod s)}}}
 \log p\cdot e(p\eta)
$$
by
$$
\frac1{\phi([q,s])}\sum_{1 \leq c \leq s \atop{(c,s)=1\atop{(q,s)|c-a}}}
e\left(\frac{bc}s\right)
\sum_{X<n\leq X+Y} e(n\eta),
$$
and the function  $S_2(\alpha)$ by
$$
\frac1{\phi(s)}
\sum_{1 \leq c \leq s \atop{(c,s)=1}}
e\left(\frac{bc}s\right)
\sum_{n\leq x}
e(n\eta)
=
\frac{\mu(s)}{\phi(s)}\sum_{n\leq x} e(n\eta).
$$
Here,  $\mu(\cdot)$ denotes  the M\"{o}bius function.

The estimates for the error terms  in our  resulting analogue of \cite[(7)]{PePi}  are identical to
the ones described in \cite[Section 4]{PePi}.
For the main term,  the only difference is in the singular series, which now originates in
\begin{eqnarray}
\sum_{s\leq\log^Cx}
\frac{\mu(s)}{\phi(s)\phi([s,q])}
\sum_{1 \leq b \leq s \atop{(b,s)=1}}
e\left(\frac{-rb}s\right)
\sum_{1 \leq c \leq s \atop{(c,s)=1\atop{(q,s)|c-a}}}
e\left(\frac{bc}s\right).\label{my26}
\end{eqnarray}
We proceed as follows.

The standard argument via the Chinese Remainder Theorem shows that the function
$$
F(s;r;q,a)
:=
\sum_{1 \leq b \leq s \atop{(b,s)=1}} e\left(\frac{-rb}s\right)
\sum_{1 \leq c \leq s \atop{(c,s)=1\atop{(q,s)|c-a}}}
e\left(\frac{bc}s\right)
$$
is multiplicative in $s$.
Indeed, for $s=uv$, $(u,v)=1$,
$u\overline u\equiv1 (\mod v)$ and $v\overline v\equiv1(\mod u)$,
note that  the relation  $c=gu\overline u+hv\overline v$ establishes a bijection between
the reduced residue classes $c$ modulo $uv$ and the pairs of reduced
residue classes $g$ modulo $v$, $h$ modulo $u$. Similarly,  the relation $b=du+fv$
establishes a bijection between the reduced residue classes $b$
modulo $uv$ and the pairs of reduced residue classes $d$ modulo $v$,
$f$ modulo $u$. Thus
$$
e\left(\frac{bc-rb}s\right)
=
e\left(\frac{(du+fv)(gu\overline u+hv\overline
v)-r(du+fv)}{uv}\right)
=
e\left(\frac{dg-rd}v\right)e\left(\frac{fh-rf}u\right).
$$
Moreover,  $(q,uv)|c-a$ if and only if $(q,u)|h-a$ and $(q,v)|g-a$.

Observe that we are only interested in square-free $s$, thus it is enough to know $F(p;r;q;a)$
for a prime $p$. A routine computation shows that
$$
F(p;r;q,a)=
\sum_{1 \leq b \leq p-1}\sum_{1 \leq c \leq p-1\atop{(q,p)|c-a}}
e\left(\frac{bc-br}p\right)=
\left\{
\begin{array}{cc}
p-1 & \mbox{ if  $p|q,\,p|a-r$,}
\\
-1 & \mbox{ if  $p|q, p\nmid a-r$,}
\\
-p+1 & \mbox{ if  $p\nmid q, p|r,$}
\\
1 & \mbox{ if  $p\nmid q, p\nmid r.$}
\end{array}
\right.
$$

Now one can easily check that,  after extending the sum over $s$ in
(\ref{my26}) up to infinity,  we have
$$
\frac1{\phi(q)}\sum_{s \geq 1}
\frac{\mu(s)}{\phi(s)} \cdot \frac{\phi(q)}{\phi([s,q])}\,F(s;r;q,a)=
\frac1{\phi(q)}\prod_{p}\left(1-\frac{\phi(q)F(p;r;q,a)}{(p-1)\phi([p,q])}\right)
=\mathfrak S(r,q,a).
$$
Proposition \ref{SiegelWalfisz} then follows. \qed


\bigskip\noindent {\bf{Proof of Theorem \ref{thmAB}.}}
Let us observe that
the expected density of twin primes  of distance $r$ is $\mathfrak S (r)$.
If  $a (\mod q)$ is an admissible residue class,
then the expected density of twin primes  of distance $r$
in the residue class $a (\mod q)$
should satisfy
$$\mathfrak S (r,q,a)=\frac{\mathfrak S(r)}{\rho(r,q)}$$
 whenever $(a,q)=(a-r,q)=1$,
where
$$\rho(r,q) := \#\{a(\mod q): (a,q)=(a-r,q)=1\}.$$

To see this, let us evaluate $\rho(r, q)$. On one hand,  we have that this
function is multiplicative in the second variable $q$. Indeed,  let
$q=uv$ with  $(u,v)=1$  and note that by the Chinese Remainder Theorem,
the relation
$b=cu\overline u+dv\overline v$ establishes a bijection between
the reduced residue classes $b$ modulo  $uv$ and the pairs of reduced
residue classes $c$ modulo  $v$, $d$ modulo  $u$. Here,
$u\overline u\equiv1(\mod v)$ and $v\overline v\equiv1 (\mod u)$. The
multiplicativity then follows from the fact that $(b-r,uv)=1$ if and
only if $(c-r,v)=(d-r,u)=1$. On the other hand, we have that
$$
\rho(r,p^\alpha)=
\ds\sum_{ 1 \leq b \leq p^{\alpha}  \atop{ p\nmid b \atop{p \nmid b-r} } } 1
=
\left\{
\begin{array}{cc}
p^\alpha-p^{\alpha-1} & \mbox{if $p|r$},
\\
p^\alpha-2p^{\alpha-1} &\mbox{ if  $p\nmid r$.}
\end{array}\right.
$$
Then our claim  follows from the equations
\begin{eqnarray}
\rho(r,q)=
\prod_{p^\alpha\|q,\,p|r}(p^\alpha-p^{\alpha-1})
\cdot
\prod_{p^\alpha\|q,p\nmid r}(p^\alpha-2p^{\alpha-1})
=
 \phi(q)
\prod_{p|q,p\nmid r}\frac{p-2}{p-1},
\label{my27}
\\
\nonumber
\frac{\mathfrak S(r)}{\rho(r,q)}
=
\frac{\mathfrak S(r)}{\phi(q)}
\prod_{p|q,p\nmid r}\frac{p-1}{p-2}
=
\frac2{\phi(q)}
\prod_{p\neq2}\frac{p(p-2)}{(p-1)^2}
\cdot
\prod_{p|r}\frac{p-1}{p-2}
\cdot
\prod_{p|q,p\nmid r}\frac{p-1}{p-2}
=
\mathfrak S(r,q,a).
\end{eqnarray}

Now let us extend the definition of $\rho(r,q)$ to characters modulo $q$:
if  $\chi$ is any non-principal character modulo $q$,  let
$$\rho(r,\chi) :=
 \ds\sum_{1 \leq b \leq q\atop{(b-r,q)=1}} \chi(b)
 =
 \ds\sum_{1 \leq b \leq q} \chi(b)\chi_0(b-r);
 $$
if $\chi_0$ denotes the principal character modulo $q$,  let
$$\rho(r,\chi_0) := \rho(r,q).$$

By the orthogonality of characters, we obtain that
\begin{eqnarray}
\mathfrak S (r,q,a)=\frac{\mathfrak S
(r)}{\rho(r,q)}=\frac{\mathfrak S
(r)}{\phi(q)}\sum_{\chi}\overline\chi(a)
\frac{\rho(r,\chi)}{\rho(r,q)},\label{my28}
\end{eqnarray}
where this formula also  incorporates all the  conditions that  $2|r$ and
$(a,q)=(a-r,q)=1$. This simple representation plays a crucial role
in the following computation.

Let $\varepsilon, M >0$, $N>N(M)$, $x>x(\varepsilon,M)$,
$x^{1/3+\varepsilon}\leq R\leq x$, $Q\leq x\log^{-N}x$ and $0\leq
X<X+Y\leq x$ be fixed, as in the statement of Theorem \ref{thmAB}. We
define, for any (even) integer $r$ and character $\chi$,
$$
F(r,\chi) :=
\ds\sum_{X<p\leq X+Y\atop{p-p'=r}} \chi(p) \log p \cdot \log p'.
$$
Then, from the orthogonality of characters we obtain that
$$
\ds\sum_{X<p\leq X+Y\atop{p\equiv a(\mod q)\atop{p-p'=r}}}
 \log p \cdot \log p'
 = \frac1{\phi(q)} \sum_{\chi} \overline \chi(a) F(r,\chi).
$$

Usually,  the main term comes from the principal character  and the
contribution of the rest is small due to the oscillation of the
characters. Unfortunately, our situation above  is more complex  and
we need to compute a dispersion over {\it{all}} characters, as follows.

The left hand
side of Theorem \ref{thmAB} can be transformed  (using
(\ref{my28}) and orthogonality) into:

\begin{eqnarray*}
S &:=&
\ds\sum_{0<r\leq R}
\ds\sum_{q\leq Q}
\ds\sum_{1 \leq a \leq q}
\left|
\ds\sum_{X<p\leq X+Y\atop{p\equiv a(\mod q)\atop{p-p'=r}}}
\log p \cdot \log p'
-
\mathfrak S(r,q,a)Y
\right|^2
\\
&=&
\ds\sum_{0<r\leq R}
\ds\sum_{q\leq Q}
\ds\sum_{1 \leq a \leq q}
\left|
\frac1{\phi(q)} \ds\sum_{\chi}
\overline \chi(a) \left( F(r,\chi)
-
\frac{\mathfrak S(r)\rho(r,\chi)Y}{\rho(r,q)}
\right)
\right|^2
\\
&=&
\ds\sum_{0<r\leq R}
\ds\sum_{q\leq Q}
\frac1{\phi(q)}
\ds\sum_{\chi}
\left|
F(r,\chi)-\frac{\mathfrak S(r)\rho(r,\chi)Y}{\rho(r,q)}
\right|^2
\\
&=&
\ds\sum_{0<r \leq R}
\ds\sum_{q  \leq Q}
\frac1{\phi(q)}
\ds\sum_{\chi}\left|F(r,\chi)
\right|^2
\\
&&-
\ds\sum_{0<r\leq R}
\ds\sum_{q\leq Q}
\frac1{\phi(q)}
\ds\sum_{\chi}
\frac{\mathfrak S(r)Y}{\rho(r,q)}2\Re\big( F(r,\chi)\rho(r,\overline\chi)\big)
\\&&+
\ds\sum_{0<r\leq R}
\ds\sum_{q\leq Q}\frac1{\phi(q)}
\ds\sum_{\chi}\frac{\mathfrak S(r)^2|\rho(r,\chi)|^2Y^2}{\rho(r,q)^2}.
\end{eqnarray*}

By the orthogonality of characters,
$$
\frac1{\phi(q)}\sum_{\chi}|\rho(r,\chi)|^2
=
\sum_{1 \leq b \leq q \atop{(b-r,q)=1}}
\sum_{1 \leq c \leq q \atop{(c-r,q)=1}}
\frac1{\phi(q)}\sum_{\chi} \chi(b)\overline\chi(c)=\rho(r,q),
$$
and  then the last term in $S$ simplifies to
$$
\sum_{0<r\leq R}\sum_{q\leq
Q}\frac1{\phi(q)}\sum_{\chi}\frac{\mathfrak S
(r)^2|\rho(r,\chi)|^2Y^2}{\rho(r,q)^2}= Y^2 \sum_{0<r\leq
R}\sum_{q\leq Q}\frac{\mathfrak S(r)^2}{\rho(r,q)}.
$$

More importantly,
\begin{eqnarray*}
&&\frac1{\phi(q)}\sum_{\chi}F(r,\chi)\rho(r,\overline\chi)\\&&=
\sum_{X<p\leq X+Y\atop{p-p'=r}}\sum_{1 \leq b \leq q \atop{(b-r,q)=1}}
\log p \cdot \log p' \frac1{\phi(q)}
\sum_{\chi}\chi(p)\overline\chi(b)\\&&=
\sum_{X<p\leq X+Y\atop{p-p'=r\atop{(pp',q)=1}}}
\log p \cdot \log p'=
\sum_{X<p\leq X+Y\atop{p-p'=r}}\log p\cdot \log p'+\O(\log^2x).
\end{eqnarray*}
The expected asymptotic for this last sum is $\mathfrak S(r)Y$,
which  is indeed true on average over $r$ from the case of
$a=q=1$ of Proposition \ref{SiegelWalfisz}.

Also,  from
$\mathfrak S(rq)\ll\mathfrak S(r)\mathfrak S(q)$ and from the fact that
$\mathfrak S(r)$ is 1 on average, that is,
$$\sum_{r\leq R}\mathfrak S(r)=\O(R)$$
and
$$\sum_{r\leq R}\mathfrak S(r)^2=\O(R),$$
 we deduce that
\begin{eqnarray*}
&&\sum_{0<r\leq R}\sum_{q\leq
Q}\frac1{\phi(q)}\sum_{\chi}\frac{\mathfrak S
(r)Y}{\rho(r,q)}2\Re\big( F(r,\chi)\rho(r,\overline\chi)\big)\\&&=
2\sum_{0<r\leq R}\sum_{q\leq Q}\frac{\mathfrak S
(r)Y}{\rho(r,q)}\big(\mathfrak S(r)Y+\O(\log^2x)\big)+
\O\left(\frac{Rx^2}{\log^Mx}\right)\\&&= 2Y^2 \sum_{0<r\leq R}\sum_{q\leq
Q}\frac{\mathfrak S(r)^2}{\rho(r,q)}+\O\left(\frac{Rx^2}{\log^Mx}\right).
\end{eqnarray*}

Putting everything together, we obtain that
$$
S=\sum_{0<r\leq R}\sum_{q\leq Q}\frac1{\phi(q)}\sum_{\chi}
\left|F(r,\chi)\right|^2 - Y^2 \sum_{0<r\leq R}\sum_{q\leq
Q}\frac{\mathfrak S(r)^2}{\rho(r,q)}+\O\left(\frac{Rx^2}{\log^Mx}\right).
$$

Note that nothing deep has happened so far beside the one
application of Proposition \ref{SiegelWalfisz}; we have utilized only
the basic properties of Dirichlet characters. Now we need to show that
the first and the second terms are asymptotically equal, that
is,  we need to be exact in our computation.

As a first step we define
\begin{eqnarray}
C(f,Q):=\sum_{q\leq Q\atop{f|q}}\frac1{\phi(q)},
\end{eqnarray}
which  satisfies
\begin{eqnarray}
C(f, Q)
\ll\frac1{\phi(f)}\log Q.\label{my29}
\end{eqnarray}
We also observe that,
if  $\chi$ mod $q$ is induced by the primitive character $\chi^*$
mod $f$, then,  due to the fact that $F(r,\chi)$ is a sum over
primes,  we have
$$F(r,\chi)=F(r,\chi^*)+\O(\log^2x).$$
By rearranging the first sum in $S$ according to primitive characters and using the above,
we then  see that
$$
S
=
\sum_{0<r\leq R}\sum_{f\leq Q}
C(f,Q)
{\sum_{\chi}}^*
\left| F(r,\chi) \right|^2 - Y^2
\sum_{0<r\leq R}
\sum_{q\leq Q}
\frac{\mathfrak S(r)^2}{\rho(r,q)}+\O\left(\frac{Rx^2}{\log^Mx}\right),
$$
where
$\sum_{\chi}^*$ is a sum over all primitive characters modulo $f$.

Now for any fixed (even) $0<r\leq R$,  we use the large sieve inequality to
estimate
$$
\sum_{Q_0<f\leq Q}C(f,Q){\sum_{\chi}}^*|F(r,\chi)|^2\ll
\left(\frac Y{Q_0}+Q\right)Y\log^3x\ll \frac{x^2}{\log^Mx},
$$
where  $Q_0:=\log^{-M-3}x$ and $N\geq M+3$.
This implies that

$$
S=\sum_{0<r\leq R}\sum_{f\leq
Q_0}C(f,Q){\sum_{\chi}}^*\left|F(r,\chi)\right|^2 - Y^2
\sum_{0<r\leq R}\sum_{q\leq Q}\frac{\mathfrak S
(r)^2}{\rho(r,q)}+\O\left(\frac{Rx^2}{\log^Mx}\right).
$$

Using the notation
\begin{eqnarray}\label{notation-psi}
\psi(X,Y;r,f,b) &:=& \sum_{{{X < p \leq X+Y} \atop {p \equiv b (\mod f)}} \atop
{p-p'=r}} \log{p} \cdot \log{p'}, \\
\label{notation-E}
E(X,Y;r,f,b) &:=& \psi(X,Y;r,f,b) - {\mathfrak{S}}(r,f,b) Y,
\end{eqnarray}
we see that
\begin{eqnarray*}
F(r,\chi)
&=&
\ds\sum_{1 \leq b \leq f}
 \chi(b)
 \psi(X,Y;r,f,b)
\\&=&
\ds\sum_{1 \leq b \leq f \atop{(b-r,f)=1}}
 \chi(b)\,\mathfrak S(r,f,b)\,Y
 +
\ds\sum_{1 \leq b \leq f}
\chi(b)E(X,Y;r,f,b)
\\&=&
\frac{\mathfrak S(r)Y}{\rho(r,f)}\rho(r,\chi) + \O\left(f\max_{(b,f)=1}|E(X,Y;r,f,b)|\right).
\end{eqnarray*}

For any small $f$ and $b$,  the sum of $|E(X,Y;r,f,b)|$ over $r$ is
sufficiently small by Proposition \ref{SiegelWalfisz}; consequently,  the same is also  true
for the sum over $f$.

Additionally, we can evaluate $\rho(r,\chi)$ for a primitive character. It
is well-known that for a primitive character $\chi$ modulo  $f$ and
for all $d|f$, $d\neq f$,
we have $\sum_{1 \leq c \leq f/d}\chi(r+cd)=0$
(see \cite[Chapter 9]{Da}). Using this,  we quickly infer that
$$
\rho(r,\chi)=\sum_{1 \leq b \leq f}
\chi(b)\sum_{d|(b-r,f)}\mu(d)
=
\sum_{d|f}\mu(d)\sum_{1 \leq b \leq f \atop{b\equiv r(\mod d)}} \chi(b)
=
\mu(f)\chi(r).
$$

Putting everything together, we arrive at the equation
$$
S=\sum_{0<r\leq R}\sum_{f\leq Q_0\atop{(r,f)=1\atop{f\text{ is
square-free}}}}C(f,Q)\frac{\mathfrak S
(r)^2Y^2}{\rho(r,f)^2}{\sum_{\chi}}^*1 - Y^2 \sum_{0<r\leq
R}\sum_{q\leq Q}\frac{\mathfrak S
(r)^2}{\rho(r,q)}+\O\left(\frac{Rx^2}{\log^Mx}\right).
$$

Let us denote the number of primitive characters modulo $f$ by
$\phi^*(f)$. We note that this is a multiplicative function for which
$\phi^*(p)=p-2$. The sum over $f$ is a quickly converging sum by
(\ref{my27}) and (\ref{my29}),  so we can drop the condition $f\leq
Q_0$ for a price already paid by the error term. Writing back the
definition of $C(f,Q)$, we obtain,  after  a little rearrangement, that
$$
S=Y^2\sum_{0<r\leq R}\mathfrak S(r)^2\sum_{q\leq Q}\frac1{\phi(q)}
\left(\sum_{f|q\atop{(r,f)=1\atop{f\text{ is square-free}}}}
\frac{\phi^*(f)}{\rho(r,f)^2} - \frac{\phi(q)}{\rho(r,q)}\right)
+ \O\left(\frac{Rx^2}{\log^Mx}\right).
$$
Finally, let us  notice that by (\ref{my27}) we actually have 0 inside the
big parantheses,  as everything is multiplicative and
$$
\sum_{f|q\atop{(r,f)=1\atop{f\text{ is square-free}}}}
\frac{\phi^*(f)}{\rho(r,f)^2}= \prod_{p|q\atop{p\nmid r}}
\left(1+\frac{\phi^*(p)}{\rho(r,p)^2}\right)=
\prod_{p|q\atop{p\nmid r}} \left(1+\frac1{p-2}\right)=
\frac{\phi(q)}{\rho(r,q)}.
$$
This completes the proof of Theorem \ref{thmAB}.
\qed



\section{Proof of Proposition \ref{mainone}}
\label{main}


\hspace*{0.5cm}
This section consists of a proof of Proposition \ref{mainone}.
This is done in two parts: an estimate of the error
term in Proposition \ref{mainone}, which  relies on Theorem \ref{thmAB},
and an estimate of the main term, which  consists
mainly in  the computation of the constant $\const$.
Note that in the course of proving Proposition \ref{mainone} we can always assume
$R\geq x^{1/3+\epsilon}$ and $Y\geq\sqrt X$,  as otherwise the
error term is an obvious upper bound for all the other terms in
(\ref{EQwantedATP}).
Note also that the term $r=1$ behaves differently,  as $p+1-r$ is always prime in this case.
However, any trivial bound (obtained by dropping the primality of $p$, or by bounding the
character by 1) shows that this term is much smaller than the
error term in Proposition \ref{mainone};  therefore it can
comfortably be excluded  from any further investigation.



Note that, using notation  (\ref{notation-E}),
Theorem \ref{thmAB} can be formulated as
$$\sum_{0<|r|\leq R}\sum_{q\leq Q}\sum_{a (\mod q)}
|E(X,Y;r,q,a)|^2\ll\frac{Rx^2}{\log^M x}, $$ whenever
$x^{1/3+\varepsilon}\leq R\leq x$,   $Q\leq x\log^{-N}x$,
and $X+Y\leq x$.

Using the same notation, as well as (\ref{notation-psi}),
we rewrite the left hand side  of  (\ref{EQwantedATP}) as
\begin{eqnarray}
&& \label{expressionx} \ds\sum_{ |r| \leq R,\,r\neq1 \atop{ f \leq
V \atop{ n \leq U } } } \frac{1}{nf} \ds\sum_{a (\mod 4n)} \left(
\frac{a}{n} \right) \psi(X,Y;r-1, nf^2, (r^2-af^2)/4)
\\
&& \label{MT} \hspace{1cm} = Y \ds\sum_{ |r| \leq R,\,r\neq1
\atop{f \leq V \atop{n \leq U } } } \frac{1}{nf} \ds\sum_{a (\mod
4n)} \left( \frac{a}{n} \right) \mathfrak{S} (r-1, nf^2,
(r^2-af^2)/4)
\\
&& \label{ET} \hspace{1cm} + \ds\sum_{ |r| \leq R,\,r\neq1
\atop{f \leq V \atop{n \leq U} } } \frac{1}{nf} \ds\sum_{a (\mod
4n)} \left( \frac{a}{n} \right) E(X,Y;r-1, nf^2, (r^2-af^2)/4).
\end{eqnarray}

\subsection{Estimate of the error term in Proposition \ref{mainone}}

\hspace*{0.5cm}
In what follows, we will show  how Theorem \ref{thmAB} allows us to
control the error term (\ref{ET}). First, using the
Cauchy-Schwarz inequality, we obtain
\begin{eqnarray}
\nonumber && \ds\sum_{{ |r| \leq R,\,r\neq1} \atop{f \leq V
\atop{n \leq U } } } \frac{1}{nf} \ds\sum_{a (\mod 4n)} \left(
\frac{a}{n} \right) E(X,Y;r-1, nf^2, (r^2-af^2)/4 )
\\
\label{afterCS} && \leq \ds\sum_{f \leq V} \frac{1}{f} \left(
\ds\sum_{ {|r| \leq R} \atop {n \leq U \atop{ a (\mod 4n)} } }
\frac{1}{n^2} \right)^{1/2} \left( \ds\sum_{{|r| \leq R,\,r\neq1}
\atop {n \leq U \atop{ a (\mod 4n)} } } E^2(X,Y; r-1, nf^2,
(r^2-af^2)/4) \right)^{1/2}.
\end{eqnarray}

The first inner sum above  is estimated trivially as
\begin{eqnarray}
\label{part1} \left( \sum_{{|r| \leq R} \atop {n \leq U \atop{ a
(\mod 4n)}} } \frac{1}{n^2} \right)^{1/2} \ll R^{1/2}
\log^{1/2}{U}.
\end{eqnarray}

For the second inner sum we observe that
\begin{eqnarray*}
\ds\sum_{ |r| \leq R,\,r\neq1 \atop{n \leq U \atop{a (\mod 4n)} }
} E^2(X,Y; r-1, nf^2, (r^2-af^2)/4) \leq \ds\sum_{ |r| \leq
R\atop r\neq1} \; \ds\sum_{q \leq 4Uf^2} \; \ds\sum_{b (\mod q)}
E^2(X,Y; r-1, q, b),
\end{eqnarray*}
as for each fixed $f,r,n$, the residue classes
$$
\left\{ b = \frac{r^2-af^2}{4}:\; a (\mod  4n) \right\}
$$
cover each residue class modulo $4n$ at most once. Then, using
Theorem \ref{thmAB}, we obtain
\begin{eqnarray}
\label{part2} \ds\sum_{|r| \leq R,\,r\neq1 \atop{n \leq U \atop{a
(\mod 4n)} } } E^2(X,Y; r-1, nf^2, (r^2-af^2)/4) &\ll& \frac{R
x^{2}}{\log^{M}{x}}
\end{eqnarray}
for any $M > 0$, provided that
\begin{eqnarray}
\label{cond_on_U_2} 4UV^2 \leq x \log^{-N}{x}. \end{eqnarray}

Using the estimates (\ref{part1}) and (\ref{part2}) in (\ref{afterCS}),
we finally obtain that
\begin{eqnarray}
\nonumber
\ds\sum_{ |r| \leq R,\,r\neq1 \atop{f \leq V \atop{n
\leq U} }} \frac{1}{nf} \ds\sum_{a (\mod 4n)} \left( \frac{a}{n}
\right) E(X,Y;r-1, nf^2, (r^2-af^2)/4)
 &\ll &
 \frac{R x \log^{1/2}{U}}{\log^{M/2}{x}} \ds\sum_{f \leq V}
\frac{1}{f} \label{boundforET}
\\
&\ll&
\frac{R x}{\log^{(M-3)/2}{x}}
\end{eqnarray}
for any $M>0$. This estimates the error term (\ref{ET}), and thus the error
term of Proposition \ref{mainone} (after renaming $M$).


\subsection{Computation of the constant in  Proposition \ref{mainone}}
\label{computeconst}


\hspace*{0.5cm}
We now treat the main term (\ref{MT}) in (\ref{expressionx}),
which is essentially a computation of the
constant $\const$ in Theorem \ref{thm1}. We first analyse the sum over $n$ and $f$ of (\ref{MT})
when $r$ is a {\it{fixed}} integer. We remark that the sum is zero if $r$ is even, thus
we can assume that $r$ is odd. We also take $r \neq 1$.

Our goal in this section is to prove:
\begin{proposition}
\label{thmfirstconstant} Let $r \neq 1$ be an odd integer. Then,
\begin{eqnarray*}
\ds\sum_{{f \leq V} \atop{n \leq U} } \frac{1}{nf} \ds\sum_{a
(\mod 4n)} \left( \frac{a}{n} \right) \mathfrak{S} (r-1, nf^2,
(r^2-af^2)/4) = C_r + \O\left(\frac1{V^2}+\frac1{\sqrt U}\right),
\end{eqnarray*}
where $C_r$ is the positive constant
\begin{eqnarray*}
C_r
&:=&
\ds\sum_{f=1}^{\infty}
\ds\sum_{n=1}^\infty
\frac{1}{nf}
\ds\sum_{a (\mod 4n)}
\left( \frac{a}{n} \right)
\mathfrak{S} (r-1, nf^2, (r^2-af^2)/4)
\\
&=& \frac{4}{3} \left( \ds\prod_{{\ell \neq 2}}
\frac{\ell^2(\ell^2-2 \ell - 2)}{(\ell-1)^3 (\ell+1)} \right)
\prod_{{\ell \mid (r-1)} \atop {\ell \neq 2}} \left(1 +
\frac{\ell+1}{\ell^2-2\ell-2}\right) \prod_{{\ell \mid r(r-2)}
\atop {\ell \neq 2}} \left(1 + \frac{1}{\ell^2-2\ell-2}\right).
\end{eqnarray*}
\end{proposition}

\noindent
{\bf{Proof of Proposition \ref{thmfirstconstant}.}}
Using the definition of ${\mathfrak S}(\cdot, \cdot, \cdot)$, we rewrite the left hand side of the
desired equation in Proposition \ref{thmfirstconstant} as
\begin{eqnarray}
\label{replacehere}
2
\left(
\ds\prod_{\ell \neq 2} \frac{\ell (\ell -2)}{(\ell-1)^2}
\right)
\ds\sum_{{f \leq V} \atop {f \;\rm{odd}}}
\ds\sum_{n \leq U}
\frac{1}{nf \phi(nf^2)}
\left( \prod_{{\ell \mid nf^2 (r-1)}\atop{\ell \neq 2}}
\frac{\ell-1}{\ell-2} \right)
\;\;c_f^r(n),
\end{eqnarray}
where
$$
c_f^r(n) := {\ds\sum_{a (\mod 4n)}}^{\prime} \leg{a}{n}
$$
and $\sum^\prime$ indicates that the sum is taken over the
invertible residues $a$ modulo $4n$ such that
\begin{eqnarray*}
&& \left( (r^2-af^2)/4, nf^2 \right) = 1
\;\;\;\;\mbox{and}\;\;\;\;
\left( (r^2-af^2)/4 - (r-1), nf^2 \right) = 1 \\
\iff && \left( r^2-af^2, 4nf^2 \right) = 4
\;\;\;\;\mbox{and}\;\;\;\; \left( r^2-af^2 - 4(r-1), 4nf^2 \right)
= 4.
\end{eqnarray*}

As $r,f$ are odd and  $(r,f)$ divides $(r^2-af^2,4nf^2)$, we must
have that $(r,f)=1$; in this case,
$$ \left( r^2-af^2, 4nf^2 \right) = 4 \iff \left( r^2-af^2, 4n \right) =
4.$$ Similarly, as $(r-2,f)$ divides $(r^2-af^2-4(r-1), nf^2)=
((r-2)^2 -af^2,nf^2)$, we must have $(r-2,f)=1$; in this case,
$$ \left( (r-2)^2-af^2, 4nf^2 \right) = 4 \iff \left( (r-2)^2-af^2, 4n \right) =
4.$$ Then
\begin{eqnarray*}
c_f^r(n) =
\left\{ \begin{array}{cl} \displaystyle {\sum_{{{a (4n)^*} \atop {(r^2-af^2,
4n)=4}}\atop{((r-2)^2-af^2, 4n)=4}} \leg{a}{n}} & \;\;\;\;\;\;\; \mbox{if $(r,f)=(r-2,f)=1$,}\\
0 & \;\;\;\;\;\;\; \mbox{otherwise,}  \end{array} \right.
\end{eqnarray*}
where $a(4n)^*$ denotes invertible residue classes $a$ modulo $4n$.

To continue the proof, we need some properties of the function $c_f^r(n)$:
\begin{lemma}
\label{computecoeffs}
Let $r \neq 1$ be an odd integer and let $f$ be a positive
odd integer such that $(r,f)=(r-2,f)=1$. Let $c_f^r(n)$ be as
defined above. The following statements  hold:
\begin{enumerate}
\item
if $n$ is odd, then
$$\displaystyle
{c_f^r(n)=
\ds\sum_{{{a (\mod n)^*} \atop {(r^2-af^2, n)=1}} \atop {((r-2)^2-af^2,n)=1}}} \leg{a}{n},
$$
where $a(\mod n)^*$ denotes invertible residue classes $a$ modulo $n$;
\item
$c_f^r(n)$ is a multiplicative function of $n$;
\item
if $\ell$ is an odd prime and $(\ell, f) =1$,
then
$${c_f^r(\ell^\alpha)\over \ell^{\alpha-1}}=\left\{\begin{array}{ll}
\ell-2 & \mbox{if $\alpha$ is even and $\ell \mid r(r-2)(r-1)$,}\\
\ell-3 & \mbox{if $\alpha$ is even and $\ell \nmid r(r-2)(r-1)$,}\\
-1 & \mbox{if $\alpha$ is odd and $\ell \mid r(r-2)(r-1)$,}\\
-2 & \mbox{if $\alpha$ is odd and $\ell \nmid r(r-2)(r-1)$;}
\end{array}\right.$$
\item
if $\ell$ is an odd prime and $\ell \mid f$ (which implies that
$(\ell,r)=(\ell, r-2)=1$ by the hypotheses on $f$), then
$$
{c_f^r(\ell^\alpha)\over \ell^{\alpha-1}}
=
\left\{\begin{array}{cc}
0  & \mbox{if $\alpha$ is odd,}\\
\ell-1 & \mbox{if $\alpha$ is even};\\
\end{array}\right.$$
\item
$\displaystyle {c_f^r(2^\alpha)\over 2^{\alpha-1}}=(-1)^\alpha$.
\end{enumerate}
\end{lemma}
\begin{proof} \begin{itemize}
\item[1.] If $n$ is odd, then \begin{eqnarray*}
c_f^r(n) &=& \ds\sum_{{{a (\mod 4n)^*, \; a \equiv 1 \mod{4}} \atop
{(r^2-af^2, n)=1}} \atop {((r-2)^2-af^2,n)=1}} \leg{a}{n} =
\ds\sum_{{{a (\mod n)^*} \atop {(r^2-af^2, n)=1}} \atop
{((r-2)^2-af^2,n)=1}} \leg{a}{n},
\end{eqnarray*}
where the last equality follows from the Chinese Remainder Theorem
and the fact that $\ds\leg{a_1}{n} = \ds\leg{a_2}{n}$ when $a_1
\equiv a_2 (\mod{n})$ for $n$ odd.
\item[2.] Let $n_1, n_2$ be two co-prime positive integers with $n_1$
odd, and let $n=n_1 n_2$. Then, using the Chinese Remainder Theorem, we obtain
\begin{eqnarray*}
c_f^r(n_1) c_f^r(n_2) &=&  \ds\sum_{{{a_1 (\mod n_1)^*} \atop
{(r^2-a_1 f^2, n_1)=1}} \atop {((r-2)^2-a_1 f^2,n_1)=1}}
\leg{a_1}{n_1} \times \ds\sum_{{{a_2 (\mod 4n_2)^*} \atop {(r^2-a_2
f^2,
4n_2)=4}} \atop {((r-2)^2-a_2 f^2,4n_2)=4}} \leg{a_2}{n_2}\\
&=&  \ds\sum_{{{a (\mod 4n_1 n_2)^*} \atop {(r^2-af^2, 4n_1 n_2)=4}}
\atop {((r-2)^2-af^2,4n_1 n_2)=4}} \leg{a}{n_1} \leg{a}{n_2}  =
c_f^r(n_1 n_2).
\end{eqnarray*}
\item[3.] We have that
\begin{eqnarray}\label{taglemma9}
c_f^r(\ell^\alpha) =  \ell^{\alpha-1} \ds\sum_{{{a (\mod \ell)^*}
\atop {(r^2-af^2, \ell)=1}} \atop {((r-2)^2-af^2,\ell)=1}}
\leg{a}{\ell}^\alpha  = \ell^{\alpha-1} \left( \ds\sum_{{{a (\mod
\ell)^*}}} \leg{a}{\ell}^\alpha - \ds\sum_{{{a (\mod \ell)^*}
\atop {a \equiv \bar{f}^{-2} r^2 (\mod{\ell)}}\text{ or}} \atop a
\equiv \bar{f}^{-2} (r-2)^2 (\mod{\ell)}} \leg{a}{\ell}^\alpha
\right),
\end{eqnarray}
where $\bar{f}$ denotes the inverse of $f$  modulo $\ell$.  We then need
to count the number of invertible residues $a$ modulo $\ell$
which are eliminated by the two congruence conditions of the second
sum. If $r \equiv 0,1,2 (\mod 4)$, there is exactly one such
residue (notice that $r^2 \equiv (r-2)^2 (\mod \ell) \iff r \equiv
1 (\mod \ell)$). In all three cases, this residue is an invertible
square modulo $\ell$,  and the second sum on the right hand side  of
(\ref{taglemma9}) has value $+1$. The result follows immediately
when  $\alpha$ is even, and follows from the orthogonality relations
when $\alpha$ is odd. If $r \not\equiv 0,1,2 (\mod 4)$, there are
two invertible residues which are eliminated by the two
congruence conditions on $a$, and the second sum on the right hand side  of
(\ref{taglemma9}) has value $+2$. The result follows as above.
\item[4.] We have that
\begin{eqnarray*}
c_f^r(\ell^\alpha) =  \ell^{\alpha-1} \ds\sum_{{{a (\mod \ell)^*}
\atop {(r^2-af^2, \ell)=1}} \atop {((r-2)^2-af^2,\ell)=1}}
\leg{a}{\ell}^\alpha = \ell^{\alpha-1} \ds\sum_{{{a (\mod \ell)^*}}}
\leg{a}{\ell}^\alpha
\end{eqnarray*}
since $(r^2-af^2, \ell)=((r-2)^2-af^2,\ell)=1$ for all $a$ when
$\ell \mid f$ and $(r,f)=(r-2,f)=1$. The result follows immediately
when $\alpha$ is even, and using the orthogonality relations when
$\alpha$ is odd.
\item[5.] Let $\alpha \geq 1$. Since $\ds\leg{a}{2}$ is a character modulo 8, we
write
\begin{eqnarray*}
c_f^r(2^\alpha) = 2^{\alpha-1} \ds\sum_{{{a (\mod 8)^*} \atop
{(r^2-af^2,2^{\alpha+2})=4}} \atop {((r-2)^2-af^2,2^{\alpha+2})=4}}
\leg{a}{2}^\alpha = 2^{\alpha-1} \leg{5}{2} = 2^{\alpha-1}
(-1)^\alpha.
\end{eqnarray*}
\end{itemize}
\end{proof}

Using parts {3.} and {4.} of Lemma \ref{computecoeffs}, we can write
\begin{eqnarray*}
&&
\ds\sum_{{f \leq V} \atop {f \;\rm{odd}}}
\ds\sum_{n \leq U}
\frac{1}{nf \phi(nf^2)}
\left(\ds\prod_{{\ell \mid nf^2 (r-1)}\atop {\ell \neq 2}} \frac{\ell-1}{\ell-2}\right)
\;\;c_f^r(n)
\\
&& =\hspace{1cm} \ds\sum_{{f=1} \atop
{(2,f)=(r,f)=(r-2,f)=1}}^\infty \ds\sum_{n=1}^{\infty}
\frac{1}{nf \phi(nf^2)} \left(\ds\prod_{{\ell \mid nf^2
(r-1)}\atop {\ell \neq 2}} \frac{\ell-1}{\ell-2}\right)
\;\;c_f^r(n) + \O\left(\frac1{V^2}+\frac1{\sqrt U}\right)
\\
&& =: D_r + \O\left(\frac1{V^2}+\frac1{\sqrt U}\right)
\end{eqnarray*}
as in \cite{DaPa}. Note that
$$C_r=2\prod_{\ell\neq2}\frac{\ell(\ell-2)}{(\ell-1)^2}\,D_r.$$

The rest of this section consists of writing $D_r$ as an Euler
product. We first write the sum over $n$ as a product. In order to
have multiplicative functions of $n$, we use the formulas
\begin{eqnarray*}
\phi(n f^2)
&=&
\frac{\phi(n) \phi(f^2) (n, f^2)}{\phi((n,f^2))},
\\
\ds\prod_{{\ell \mid nf^2 (r-1)}\atop {\ell \neq 2}}
\frac{\ell-1}{\ell-2}
&=&
\left(
\ds\prod_{{\ell \mid n}\atop {\ell \neq 2}}
\frac{\ell-1}{\ell-2}
\right)
\left(
\ds\prod_{{\ell \mid f^2 (r-1)}\atop {\ell\neq 2}}
\frac{\ell-1}{\ell-2}
\right)
\left(
\ds\prod_{{\ell \mid (n,f^2 (r-1))}\atop {\ell\neq 2}}
\frac{\ell-2}{\ell-1}
\right).
\end{eqnarray*}
Now we rewrite $D_r$ as
\begin{eqnarray*}
\left(
\ds\sum_{{{f=1} \atop {f \;\rm{odd}}} \atop {(r,f) = (r-2, f)=1}}^\infty
\frac{1}{f \phi(f^2)}
\ds\prod_{{\ell \mid f^2 (r-1)}\atop {\ell \neq 2}}
\frac{\ell-1}{\ell-2}
\right)
\ds\sum_{n=1}^\infty
\frac{c_f^r(n)}{n \phi(n)} \frac{\phi((f^2,n))}{(f^2,n)}
\left(
\ds\prod_{{\ell \mid n}\atop {\ell \neq 2}}
\frac{\ell-1}{\ell-2}
\right)
\left(
\ds\prod_{{\ell \mid (n,f^2(r-1))}\atop {\ell \neq 2}}
\frac{\ell-2}{\ell-1}
\right).
\end{eqnarray*}

The sum over $n$ is the sum of a multiplicative function of $n$
whose factors at prime powers $\ell^{\alpha}$ depend on the
divisibilities of $f, r, r-1$ and $r-2$ by $\ell$. Using Lemma
\ref{computecoeffs}, we can then write the $n$-sum as
\begin{eqnarray*}
\left(
\ds\prod_{\ell \mid f}
\ds\sum_{\alpha=0}^{\infty}
a_r(\ell^\alpha)
\right)
\left(
\ds\prod_{\ell \nmid f}
\ds\sum_{\alpha=0}^{\infty}
b_r(\ell^\alpha)
\right)
&=&
\left(
\ds\prod_{\ell}
\ds\sum_{\alpha=0}^{\infty}
b_r(\ell^\alpha)
\right)
\left(
\ds\prod_{\ell \mid f}
\frac{\sum_{\alpha=0}^{\infty} a_r(\ell^\alpha)}{\ds\sum_{\alpha=0}^{\infty} b_r(\ell^\alpha)}
\right),
\end{eqnarray*}
where $a_r(1)=b_r(1)=1$ and for $\ell \neq 2$ and $\alpha \geq 1$,
\begin{eqnarray*}
a_r(\ell^\alpha)
&=&
\left\{ \begin{array}{cc}
0 & \mbox{if $\alpha$ odd},
\\
(\ell-1)/(\ell^{\alpha+1}) & \mbox{if $\alpha$ even};
\end{array} \right.
\\
b_r(\ell^\alpha)
&=&
\left\{ \begin{array}{cc}
-2 / \ell^{\alpha}(\ell-2) & \mbox{if $\alpha$ odd and $\ell \nmid r(r-1)(r-2)$},
\\
(\ell-3) / \ell^\alpha (\ell-2) & \mbox{if $\alpha$ even and $\ell
\nmid r(r-1)(r-2)$},
\\
-1/\ell^\alpha(\ell-2) & \mbox{if $\alpha$ odd and $\ell \mid r(r-2)$},
\\
1/\ell^\alpha & \mbox{if $\alpha$ even and $\ell \mid r(r-2)$},
\\
-1/\ell^\alpha(\ell-1) & \mbox{if $\alpha$ odd and $\ell \mid r-1$},\\
(\ell-2)/\ell^\alpha(\ell-1) & \mbox{if $\alpha$ even and $\ell \mid
r-1$}.
\end{array} \right.
\end{eqnarray*}

Replacing in $D_r$, this gives
\begin{eqnarray*}
D_r &=& \left( \ds\prod_{\ell} \sum_{\alpha=0}^{\infty}
b_r(\ell^\alpha) \right) \ds\sum_{{{f=1} \atop {f \;\rm{odd}}} \atop
{(r,f) = (r-2, f)=1}}^\infty \frac{1}{f \phi(f^2)} \left(
\ds\prod_{{\ell \mid f^2 (r-1)}\atop{\ell \neq 2}}
\frac{\ell-1}{\ell-2} \right) \left( \ds\prod_{\ell \mid f}
\frac{\ds\sum_{\alpha=0}^{\infty}
a_r(\ell^\alpha)}{\ds\sum_{\alpha=0}^{\infty} b_r(\ell^\alpha)}
\right)
\\
&=& \left( \ds\prod_{\ell} \ds\sum_{\alpha=0}^{\infty}
b_r(\ell^\alpha) \right) \left( \ds\prod_{{\ell \mid
(r-1)}\atop{\ell \neq 2}} \frac{\ell-1}{\ell-2} \right)\times \\
&& \hspace{.1in} \times \ds\sum_{{{f = 1} \atop {f \;\rm{odd}}}
\atop {(r,f) = (r-2, f)=1}}^\infty \frac{1}{f \phi(f^2)} \left(
\ds\prod_{{\ell \mid f^2}\atop{\ell \neq 2}} \frac{\ell-1}{\ell-2}
\right) \left( \ds\prod_{{\ell \mid(f^2,r-1)}\atop{\ell \neq 2}}
\frac{\ell-2}{\ell-1} \right) \left( \ds\prod_{\ell \mid f}
\frac{\ds\sum_{\alpha=0}^{\infty}
a_r(\ell^\alpha)}{\ds\sum_{\alpha=0}^{\infty} b_r(\ell^\alpha)}
\right).
\end{eqnarray*}
The sum over $f$ in the last expression is a sum of multiplicative functions of $f$, which  we
write  as
$$
\ds\prod_{\ell \nmid 2r(r-2)}
\ds\sum_{\alpha=0}^{\infty} c_r(\ell^\alpha).
$$
Here $c_r(1)=1$ and for any prime $\ell$
with $\ell \nmid 2r(r-2)$ and any $\alpha \geq 1$, we have
\begin{eqnarray*}
c_r(\ell^\alpha) &=& \left\{ \begin{array}{ll}
\displaystyle{\frac{1}{\ell^{3\alpha-1} (\ell-2)}
\cdot
\frac{\ds\sum_{\beta \geq 0} a_r(\ell^\beta)}{\ds\sum_{\beta \geq 0} b_r(\ell^\beta)}}
& \mbox{if $\ell \nmid r-1$,}\\
\\
\displaystyle{\frac{1}{\ell^{3\alpha-1} (\ell-1)}
\cdot
\frac{\ds\sum_{\beta \geq 0} a_r(\ell^\beta)}{\ds\sum_{\beta \geq 0} b_r(\ell^\beta)}}
& \mbox{if $\ell \mid r-1$.}
\end{array} \right.
\end{eqnarray*}
Then,

\begin{eqnarray}
\label{C1}  D_r &=& \left( \ds\prod_{\ell}
\ds\sum_{\alpha=0}^{\infty} b_r(\ell^\alpha) \right) \left(
\ds\prod_{{\ell \mid (r-1)}\atop{\ell\neq2}} \frac{\ell-1}{\ell-2}
\right) \left( \ds\prod_{\ell \nmid 2r(r-2)}
\ds\sum_{\alpha=0}^{\infty} c_r(\ell^\alpha) \right).
\end{eqnarray}

We now compute the sums appearing in (\ref{C1}), using the formulas
for $a_r(\ell)$ and $b_r(\ell)$ listed above:
\begin{lemma}
\label{formulas}
Let $\ell$ be an odd prime and $\alpha \geq 1$. Let
$a_r(\ell^\alpha)$, $b_r(\ell^\alpha)$ and $c_r(\ell^\alpha)$ be as defined
above. We have:
\begin{enumerate}
\item
$
A(\ell)
:=
\displaystyle{\sum_{\alpha=0}^{\infty} a_r(\ell^\alpha)
= \frac{\ell^2 + \ell + 1}{\ell(\ell+1)}}$;
\item
if $\ell \nmid r(r-1)(r-2)$, then
$
B^{(1)}(\ell) :=
\displaystyle{\sum_{\alpha=0}^{\infty} b_r(\ell^\alpha) = \frac{\ell^3
- 2\ell^2 - 2\ell - 1}{(\ell-2)(\ell^2-1)}}$;
\item
if $\ell \mid r-1$, then $ B^{(2)}(\ell) :=
\displaystyle{\sum_{\alpha=0}^{\infty} b_r(\ell^\alpha) =
\frac{\ell^3 - \ell^2 - \ell - 1}{(\ell-1)^2(\ell+1)}}$;
\item
if $\ell \nmid r-1$, but $\ell \mid r(r-2)$, then
$
B^{(3)}(\ell) :=
\displaystyle{\sum_{\alpha=0}^{\infty} b_r(\ell^\alpha) =
\frac{\ell(\ell^2- 2\ell - 1)}{(\ell-2)(\ell^2-1)}}$;
\item
if $\ell \nmid 2r(r-2)(r-1)$, then
$
C^{(1)}(\ell) :=
\displaystyle{\sum_{\alpha=0}^{\infty}} c_r(\ell^\alpha) =
\frac{\ell^3-2\ell^2-2\ell}{\ell^3-2\ell^2-2\ell-1}$;
\item
if $\ell \mid r-1$, then $ C^{(2)}(\ell) :=
\displaystyle{\sum_{\alpha=0}^{\infty}} c_r(\ell^\alpha) =
\frac{\ell(\ell^2-\ell-1}{\ell^3-\ell^2-\ell-1}$;
\item
if $\ell=2$, then
$B(2) := \displaystyle{\sum_{\alpha=0}^{\infty}} b_r(\ell^\alpha) = \frac{2}{3}.$
\end{enumerate}
\end{lemma}
\begin{proof} All the computations are straightforward, following in
one line from the formula for the sum of the geometric series.
\end{proof}

We extend the definitions of
$A(\ell), B^{(1)}(\ell), B^{(2)}(\ell), B^{(3)}(\ell), C^{(1)}(\ell), C^{(2)}(\ell)$
introduced in  Lemma \ref{formulas} to any odd prime $\ell$
(independently of the relation between $\ell$ and $r$).
Then we rewrite $D_r$ as
\begin{eqnarray*}
D_r
&=&
B(2)
\left(
\ds\prod_{\ell \neq 2}
\ds\sum_{\alpha \geq 0}
b_r(\ell^\alpha)
\right)
\left(
\ds\prod_{{\ell \mid r-1}\atop{\ell\neq2}}
\frac{\ell-1}{\ell-2}
\right)
\left(
\ds\prod_{\ell \nmid 2r(r-2)}
\ds\sum_{\alpha \geq 0}
c_r(\ell^\alpha)
\right)
\\
&=&
B(2)
\left(
\ds\prod_{\ell \neq 2} B^{(1)}(\ell) C^{(1)}(\ell)
\right)
\left(
\ds\prod_{{\ell \mid r(r-2)}\atop{\ell \neq 2}}
B^{(1)}(\ell)^{-1} C^{(1)}(\ell)^{-1} B^{(3)}(\ell)
\right)
\times
\\
&&
\hspace{1cm} \times
\left(
\ds\prod_{{\ell \mid r-1}\atop{\ell \neq 2}}
B^{(1)}(\ell)^{-1} C^{(1)}(\ell)^{-1} B^{(2)}(\ell) \frac{\ell-1}{\ell-2} C^{(2)}(\ell)
\right).
\end{eqnarray*}
Using Lemma \ref{formulas}, we compute
\begin{eqnarray*}
B^{(1)}(\ell) C^{(1)}(\ell) &=& \frac{\ell^3-2\ell^2-2 \ell}{\ell^3-2\ell^2-\ell+2} =
\frac{\ell(\ell^2-2\ell-2)}{(\ell-2)(\ell^2-1)},\\
B^{(1)}(\ell)^{-1} C^{(1)}(\ell)^{-1} B^{(3)}(\ell) &=& \frac{\ell^2-2\ell-1}{\ell^2-2\ell-2} =
1 + \frac{1}{\ell^2-2\ell-2},\\
B^{(1)}(\ell)^{-1} C^{(1)}(\ell)^{-1} B^{(2)}(\ell) C^{(2)}(\ell)
\frac{\ell-1}{\ell-2} &=& \frac{\ell^2-\ell-1}{\ell^2-2\ell-2} = 1 +
\frac{\ell+1}{\ell^2-2\ell-2}.
\end{eqnarray*}
Finally, replacing all the above in (\ref{replacehere}), we obtain
\begin{eqnarray*}
C_r
&=&
\frac{4}{3}
\left(
\ds\prod_{{\ell \neq 2}} \frac{\ell(\ell-2)}{(\ell-1)^2}
\cdot
\frac{\ell(\ell^2-2\ell-2)}{(\ell-2)(\ell^2-1)}
\right)
\cdot
\ds\prod_{{\ell \mid r-1} \atop {\ell \neq 2}}
 \left(1 + \frac{\ell+1}{\ell^2-2\ell-2}\right)
\cdot
\ds\prod_{{\ell \mid r(r-2)} \atop {\ell \neq 2}}
\left(1 + \frac{1}{\ell^2-2\ell-2}\right)
\\
&=& \frac{4}{3} \ds\prod_{{\ell \neq 2}} \frac{\ell^2
(\ell^2-2\ell-2)}{(\ell-1)^3 (\ell+1)} \cdot \ds\prod_{{\ell \mid
(r-1)} \atop {\ell \neq 2}} \left(1 +
\frac{\ell+1}{\ell^2-2\ell-2}\right) \cdot \ds\prod_{{\ell \mid
r(r-2)} \atop {\ell \neq 2}} \left(1 +
\frac{1}{\ell^2-2\ell-2}\right).
\end{eqnarray*}
This completes the proof of Proposition \ref{thmfirstconstant}.
\qed

\subsection{The average constant}
\label{constantgalla}

\hspace*{0.5cm} Using Proposition \ref{thmfirstconstant}, the main
term (\ref{MT}) of (\ref{expressionx}) is $ Y \sum_{ |r| \leq R,\,
r\neq1 \; \text{odd}} C_r$;  thus now  we need to average the constant
$C_r$ over $r$. This calculation has similarities with the one done
by Gallagher in \cite{Ga} for the average of the standard twin
prime constant.  As such, we will follow  the notation used in \cite{Ga}.

The main result of this section is:
\begin{lemma}
As $R \rightarrow \infty$,
\label{likegallagher} $$\sum_{ |r| \leq R \atop{r\neq1 \;
\text{odd}}} C_r = {\mathfrak{C}} R + \O \left( \log^2{R}
\right).$$
\end{lemma}

\begin{proof}
Let us write
\begin{eqnarray}
\label{replaceattheend} C_r = \frac{4}{3} \ds\prod_{\ell \neq 2}
\frac{\ell^2 (\ell^2-2\ell-2)}{(\ell-1)^3 (\ell+1)} \cdot \ds\prod_{
\ell \neq 2 \atop{ \ell | r (r-1) (r-2) } } \left( 1 + e_r(\ell)
\right),
\end{eqnarray}
where
\begin{eqnarray*}
e_r(\ell)
&:=&
\left\{ \begin{array}{cc}
e^{(1)}(\ell) & \mbox{if $\ell \mid r-1$} \\
e^{(2)}(\ell) & \mbox{if $\ell \mid r(r-2)$}\\
0 & \mbox{otherwise}
\end{array} \right.
=
\left\{ \begin{array}{cc}
\ds\frac{\ell+1}{\ell^2-\ell-2} & \mbox{if $\ell \mid r-1$,} \\
\ds\frac{1}{\ell^2-2\ell-2} & \mbox{if $\ell \mid r(r-2)$,}\\
0 & \mbox{otherwise.}
\end{array} \right.
\end{eqnarray*}
Let us also fix the following notation: for $r\neq1$ odd,  we take
\begin{eqnarray*}
{\cal P}(r)
&:=&
\left\{ \ell \;\mbox{odd prime}:\; \ell \mid r(r-1)(r-2) \right\},
\\
{\cal F}(r)
&:=&
\left\{ q \;\mbox{positive square-free integer}:\; \ell \mid q \Rightarrow
\ell \in {\cal P}(r) \right\},
\\
{\cal D}(R) &:=& \cup_{ |r| \leq R \atop{r\neq1 \; \text{odd}} }
{\cal F}(r).
\end{eqnarray*}
We want to evaluate
\begin{eqnarray}
\label{defineS}
{\cal{S}}
 := \ds \sum_{|r| \leq R \atop{r\neq1 \;
\text{odd}} } \prod_{\ell \neq 2 \atop{\ell | r(r-1)(r-2)} } (1 +
e_{r}(\ell))=  \ds\sum_{|r| \leq R \atop{r\neq1 \; \text{odd}} }
\ds\sum_{q \in {\cal F}(r)} e_r(q),
\end{eqnarray}
where $e_r(1) = 1$ and,  for $q \neq 1$, $e_r(q) = \ds\prod_{\ell|q}
e_{r}(\ell)$. We write
\begin{eqnarray*}
{\cal{S}}
&=& \ds\sum_{q \in {\cal D}(R)} \ds\sum_{ |r| \leq R \atop{
r\neq1 \; \text{odd}} } e_r(q) = \ds\sum_{q \in {\cal D}(R)}
\ds\sum_{\text{all possible} \atop{e=e(q)} } \ds\sum_{ |r| \leq R
\atop{ r\neq1 \; \text{odd} \atop{e_r(q) = e} } } e_r(q)
\\
&=& \ds\sum_{q \in {\cal D}(R)} \ds\sum_{\text{all possible}
\atop{e=e(q)} } \#\{ |r| \leq R: \; r\neq1 \; \text{odd}, e_r(q)
= e \}
\\
&=& \ds\sum_{q \in {\cal D}(R)} \ds\sum_{v=v(q)} \ds\prod_{\ell |
q} e^{v(\ell)}(\ell) N(q, v),
\end{eqnarray*}
where the sum $\ds\sum_{v=v(q)}$ is over all maps
$v: \{\ell: \; \ell | q\} \longrightarrow \{1, 2\}$ and where
$$
N(q, v) := \# \{ |r| \leq R:\; r\neq1 \; \text{odd}, e_r(\ell) =
e^{v(\ell)}(\ell) \; \forall \ell | q \}.
$$

By looking at the conditions imposed on $\ell$ when defining
$e^{(1)}(\ell)$ and $e^{(2)}(\ell)$, we see that $N(q, v)$ is the
number of integers $|r| \leq R$ with $r\neq1 \;\text{odd}$ such
that
\begin{eqnarray*}
r &\equiv& 1 (\mod 2),
\\
r &\equiv& 1 (\mod \ell) \; \forall \ell | q \; \text{with $v(\ell) = 1$},
\\
r &\equiv& 0 \; \text{or} \; 2 (\mod \ell) \; \forall \ell | q \;
\text{with $v(\ell) = 2$}.
\end{eqnarray*}
Therefore, by using the Chinese Remainder Theorem, $r$ as above lies in one of
$\ds\prod_{\ell | q} 2^{v(\ell) - 1}$ distinct residue classes modulo $2 q$.
Consequently,
\begin{eqnarray*}
N(q, v)
&=&
\left(\ds\prod_{\ell | q} 2^{v(\ell) -1}\right)
\left(\frac{2R + 1}{2 q} + \O(1)\right)
\\
&=&
\frac{R}{q} \ds\prod_{\ell | q} 2^{v(\ell) - 1} + \O\left(2^{\omega(q)}\right),
\end{eqnarray*}
where $\omega(q)$ denotes the number of distinct prime factors of $q$.
We plug this in the formula for ${\cal{S}}$ and obtain
\begin{eqnarray*}
{\cal{S}}
&=&
R
\ds\sum_{q \in {\cal D}(R)}
\frac{1}{q}
\ds\sum_{v = v(q)}
\ds\prod_{\ell | q} e^{v(\ell)}(\ell) 2^{v(\ell)-1}
+
\O\left(
\ds\sum_{q \in {\cal D}(R)}
2^{\omega(q)}
\ds\sum_{v=v(q)}
\ds\prod_{\ell | q} e^{v(\ell)}(\ell)
\right)
\\
&=:&
{\cal{S}}_{\text{main}} + {\cal{S}}_{\text{error}}.
\end{eqnarray*}

To estimate ${\cal{S}}_{\text{main}}$, we observe that we have
$$
{\cal{S}}_{\text{main}} = R \ds\sum_{q \in {\cal D}(R)} G(q)
$$
for some multiplicative function $G(q)$. Therefore
\begin{eqnarray*}
{\cal{S}}_{\text{main}} &=& R \ds\prod_{\ell \leq R \atop{\ell \neq 2}} (1 +
G(\ell)) = R \ds\prod_{\ell \leq R \atop{\ell \neq 2}} \left(1 +
\frac{e^{(1)}(\ell) + 2 e^{(2)}(\ell)}{\ell}\right)
\\
&=& R \ds\prod_{\ell \leq R \atop{\ell \neq 2}} \frac{\ell^3 - 2
\ell^2 - \ell + 3}{\ell (\ell^2 - 2 \ell - 2)} = R \ds\prod_{\ell
\neq 2} \frac{\ell^3 - 2 \ell^2 - \ell  + 3}{\ell (\ell^2 - 2 \ell -
2)} + \O(1).
\end{eqnarray*}

Now let us estimate ${\cal{S}}_{\text{error}}$.
As for ${\cal{S}}_{\text{main}}$, we observe that we have
$$
{\cal{S}}_{\text{error}} = \O\left(\ds\sum_{q \in {\cal D}(R)} F(q)\right)
$$
for some multiplicative function $F(q)$. Therefore
\begin{eqnarray*}
{\cal{S}}_{\text{error}}
&=&
\O\left(
\ds\prod_{\ell \leq R \atop{\ell \neq 2}}
[1 + 2 (e^{(1)}(\ell) + e^{(2)}(\ell))]
\right)
\\
&=& \O\left( \ds\prod_{\ell \leq R \atop{\ell \neq 2}} \left(1 +
\frac{2(\ell^2 + \ell - 1)}{\ell(\ell^2 - 2 \ell - 2)}\right)
\right) =  \O\left((\log R)^2\right).
\end{eqnarray*}

We put the two estimates together and obtain
$$
{\cal{S}} = R \ds\prod_{\ell \neq 2} \left[
\frac{\ell^3-2\ell^2-\ell+3}{\ell(\ell^2-2\ell-2)} \right] +
\O\left((\log R)^2\right).
$$
By replacing (\ref{defineS}) in (\ref{replaceattheend}), this gives
\begin{eqnarray*}
\ds\sum_{ |r| \leq R \atop{r\neq1 \; \text{odd}} } C_r &=&
\frac{4R}{3} \ds\prod_{{\ell \neq 2}} \frac{\ell^2
(\ell^2-2\ell-2)}{(\ell-1)^3 (\ell+1)} \cdot
\frac{\ell^3-2\ell^2-\ell+3}{\ell(\ell^2-2\ell-2)} +
\O\left(\log^2R\right)\\
&=&
\frac{4R}{3}
 \ds\prod_{\ell \neq 2}
\frac{\ell^4 - 2 \ell^3 - \ell^2 + 3 \ell}{ (\ell - 1)^3 (\ell+1)} +
\O\left(\log^2{R}\right),
\end{eqnarray*}
which  completes the proof of Lemma \ref{likegallagher}.
\end{proof}

Replacing Proposition \ref{thmfirstconstant} and Lemma
\ref{likegallagher} in (\ref{MT}), this concludes the proof of
Proposition \ref{mainone}.

\end{document}